\documentclass[draft]{amsart}

\usepackage[all]{xy}
\usepackage{amssymb,amsmath}

\usepackage{amsfonts}
\usepackage{amscd}

\usepackage{url}

\makeatletter
\def\url@leostyle{%
  \@ifundefined{selectfont}{\def\UrlFont{\sf}}{\def\UrlFont{\small\ttfamily}}}
\makeatother
\newtheorem{thm}{Theorem}[section]
\newtheorem{lem}[thm]{Lemma}
\newtheorem{prop}[thm]{Proposition}
\newtheorem{cor}[thm]{Corollary}

\newtheorem{defn}[thm]{Definition}

\newtheorem{nrmk}[thm]{Remark}
\newtheorem{expl}[thm]{Example}

\newcommand{\pf}{{\bf Proof. }}

\newcommand{\into}{\longrightarrow}

\renewcommand{\tilde}{\widetilde}
\renewcommand{\bar}{\overline}

\newcommand{\NN}{\mathbb{N}}
\newcommand{\ZZ}{\mathbb{Z}}

\newcommand{\M}{\mbox{${\mathcal M}$}}

\newcommand{\G}{\mbox{${\mathcal G}$}}
\newcommand{\F}{\mbox{${\mathcal F}$}}
\newcommand{\I}{\mbox{${\mathcal I}$}}

\renewcommand{\mod}{\mathrm{Mod}}

\newcommand{\supp}{\mathrm{supp}}

\newcommand{\lind}[1]{\underset{#1}{\underrightarrow{\lim}}}
\newcommand{\Lind}{\underrightarrow{\lim}}  %grazie Anna

\begin{document}

\title {Poincar\'e - Verdier duality in o-minimal structures}

\author {M\'{a}rio J. Edmundo}

\address{ CMAF Universidade de Lisboa\\
Av. Prof. Gama Pinto 2\\
1649-003 Lisboa, Portugal}

\email{edmundo@cii.fc.ul.pt}

\author{Luca Prelli %\thanks{%Supported by NSERC. \newline
%{\it MSC}: 03C64; 55N30.
%{\it Keywords and phrases:} O-minimal structures, sheaf cohomology. }
}

\address{ %CAUL Universidade de Lisboa\\
%Av. Prof. Gama Pinto 2\\
%1649-003 Lisboa, Portugal\\
%and\\
Universit\`a di Padova\\
Dipartimento di Matematica Pura ed Applicata\\
Via Trieste 63\\
35121 Padova, Italy}
\email{lprelli@math.unipd.it}

\date{\today}
%\subjclass{Primary ; Secondary }
\thanks{The first author was supported by Funda\c{c}\~ao para a Ci\^encia e a Tecnologia, Financiamento Base 2008 - ISFL/1/209. This work was developped within project POCTI-ISFL:1-143 of CAUL supported
by FCT and FEDER.\newline
 }


\maketitle
\begin{abstract}
Here we prove a  Poincar\'e - Verdier  duality theorem for the o-minimal sheaf cohomology with definably compact supports of definably normal, definably locally compact spaces  in an arbitrary  o-minimal structure.
%\newline
%
%\noindent
%R{\tiny \'ESUM\'E.}  On d\'emontre une dualit\'e de Poincar\'e - Verdier dans le cadre de la cohomologie o-minimale des faisceaux avec support compact et d\'efinissable  sur un espace d\'efinissablement normal, d\'efinissablement localement compact dans
%associ\'e \`a une structure o-minimale arbitraire.
\end{abstract}

\newpage

\begin{section}{Introduction}\label{section introduction}
We fix an arbitrary o-minimal structure ${\mathcal M}=(M,<, \ldots )$ and  work in the category of definable spaces, $X$, in ${\mathcal M}$ with the o-minimal site on $X$, with morphisms being definable continuous maps. The o-minimal site on $X$ is the site whose
underlying category is the set of all relatively open definable subsets of $X$ (open in the strong, o-minimal topology) with morphisms the inclusions and admissible coverings being  covers by open definable sets with finite subcoverings.

The o-minimal setting generalizes the semi-algebraic and globally sub-analytic contexts (\cite{vdd}), and so our first main theorem (on Subsection \ref{section shsptop}) generalizes the existence of sheaf cohomology with supports in semi-algebraic geometry, as described in the book~\cite{D3}. This o-minimal sheaf cohomology with supports satisfies the Eilenberg-Steenrod axioms adapted to the o-minimal site -  for the homotopy axiom we need to assume that $\M$ has definable Skolem functions and the definable space $X$ involved is definably normal such that for every closed interval $[a,b]\subseteq M$ the projection $X\times [a,b]\into X$ maps closed definable subsets into closed definable subsets. Other cohomology theories have been constructed for o-minimal structures of special types in the past. Simplicial and singular cohomologies were constructed  in o-minimal expansions of fields by A.Woerheide in his doctoral thesis, a report of which can be found in~\cite{ew}. A sheaf cohomology without supports has been constructed in~\cite{ejp} for o-minimal structures (with the extra technical assumptions for the homotopy axiom given above), which generalised the sheaf cohomology without supports for real algebraic geometry of Delfs, for which he proved the homotopy axiom in~\cite{D2}. The theory presented here generalises all of these and is an extension of the corresponding theory in topological spaces (\cite{b}, \cite{g}, \cite{i} and \cite{ks1}).

Following the classical proof of the Poincar\'e - Verdier  duality for topological spaces
%presented in the book by Iversen (\cite{i})
we prove here a  version Verdier  duality theorem for the o-minimal sheaf cohomology with definably compact supports of definably normal, definably locally compact spaces  in an arbitrary  o-minimal structure (Theorem \ref{thm vd dual}).  This result is  new even in the semi-algebraic context. We do not develop yet the full theory of proper direct image and its dual in the o-minimal context but nevertheless we prove our version of Verdier duality  by considering  inclusions of definably locally closed definable subsets. The theory of proper direct image is partially developed in the semi-algebraic case in the book by Delf's (\cite{D3}). In the sub-analytic context there are several approaches to this theory by Kashiwara and Schapira (\cite{ks2}) and also L. Prelli (\cite{lucap}).

From Verdier duality we derive the Poincar\'e and Alexander duality theorems (Theorems \ref{thm poinc dual} and  \ref{thm alex dual}). The later results are based on a general and new orientation theory for definable manifolds which we show to be the same as the orientation theory in o-minimal expansions of fields defined in \cite{bo} and \cite{beo} using o-minimal singular homology. (See subsection \ref{subsection duality in fields}).

Our Poincar\'e - Verdier  duality theory relays heavily on the theory of normal and constructible supports and o-minimal cohomological $\Phi$-dimension. This rather technical theory in presented in Section \ref{section appendix} and is the o-minimal version of the topological theory of paracompactifying families of  supports and cohomological $\Phi $-dimension and generalizes the corresponding theory in the semi-algebraic context (\cite{D3}).

The motivation for developing this general o-minimal Poincar\'e - Verdier  duality is to be able to apply it to compute the o-minimal cohomology of definably compact definable groups defined in arbitrary o-minimal structures generalizing in this way  the computation of the o-minimal singular cohomology of definable groups in o-minimal expansions of fields already presented in \cite{eo}. We hope do this in a different paper.

\medskip
\textbf{ Acknowledgement: }  The contents of Subsection \ref{section shsptop} and most of Subsection \ref{subsection normal and constructible supports} were previously worked out jointly with Gareth O. Jones and Nicholas J. Peatfield. The first author  would like to thank them for allowing him to include this material here.
\end{section}

\begin{section}{Notations and review}\label{section notations}

In this section we recall some preliminaries notions about sheaves on topological spaces and the previous results about sheaves on the o-minimal spectrum of a definable space. For further details about classical sheaf theory, see for example \cite{b}, \cite{g}, \cite{i}, \cite{ks1} and \cite{ks3}.  Good references on o-minimality are, for example, the book \cite{vdd} by van den Dries and the notes \cite{c} by Michel Coste. For semi-algebraic geometry relevant to this paper the reader should consult  the work by Delfs (\cite{D2} and \cite{D3}), Delfs and Knebusch (\cite{dk2}) and the book \cite{BCR} by Bochnak, Coste and M-F. Roy.

\subsection{Sheaves on topological spaces}\label{section shtop}

Let $X$ be a topological space and let $k$ be a field. As usual, we will set $\mathrm{Mod}(k_X)$ the category of sheaves of $k$-modules on $X$. This is a Grothendieck category, hence it admits enough injectives and a family of generators (the sheaves $k_U$ defined below). Moreover filtrant inductive limits are exact.

Let $f:X \to Y$  be a morphism of topological spaces. As usual we denote by $f_*$ and $f^{-1}$ the functors of direct and inverse image. In particular, when $Y$ is a subset of $X$ we will denote by $i_Y:Y \hookrightarrow X $ the inclusion. \\

When $S$ is closed and $\mathcal{F}\in\mathrm{Mod}(k_X)$ one sets $\mathcal{F}_S=i_{S*}i_S^{-1}\mathcal{F}$, and when $U$ is open $\mathcal{F}_U=\ker(\mathcal{F}\to\mathcal{F}_{X\setminus U})$ (or equivalently $\mathcal{F}_U$ is the sheaf associated to the presheaf $V\mapsto \Gamma(V;\mathcal{F}_U)$ which is $\Gamma(V;\mathcal{F})$ if $V \subseteq U$ and $0$ otherwise). When $Z=U \cap S$ set $\mathcal{F}_Z=(\mathcal{F}_U)_S$.  The functor $(\bullet )_Z$ is exact and ${\mathcal F}_Z$ is characterized by ${\mathcal F}_{Z|Z}={\mathcal F}_{|Z}$ and ${\mathcal F}_{Z|X\setminus Z}=0$. If $Z'$ is another locally closed subset of $X$, then $({\mathcal F}_Z)_{Z'}={\mathcal F}_{Z\cap Z'}.$ When $\mathcal{F}=k_X$ is the constant sheaf on $X$ we just set $k_Z$ instead of $(k_X)_Z$. If $Z_1,Z_2$ are locally closed and $Z_1$ is closed in $Z_2$ we have an exact sequence
$$
0 \to \mathcal{F}_{Z_2 \setminus Z_1} \to \mathcal{F}_{Z_2} \to \mathcal{F}_{Z_2 \cap Z_1} \to 0.
$$

When $U$ is open one sets $\Gamma_U\mathcal{F}=i_{U*}i_U^{-1}\mathcal{F}$. Then we have $\Gamma(V;\Gamma_U\mathcal{F})=\Gamma(U \cap V;\mathcal{F})$.  When $S$ is closed $\Gamma_S\mathcal{F}=\ker(\mathcal{F}\to\Gamma_{X\setminus S}\mathcal{F})$ (sections with support in $S$). When $Z=U \cap S$ we  set $\Gamma_Z=\Gamma_U \circ \Gamma_S$. The functor $\Gamma _Z(\bullet )$ is left exact and if $Z'$ is another locally closed subset, then $\Gamma _{Z'}(\Gamma _Z{\mathcal F})=\Gamma _{Z\cap Z'}{\mathcal F}.$ If $Z_1,Z_2$ are locally closed and $Z_1$ is closed in $Z_2$ we have an exact sequence
$$
0 \to \Gamma_{Z_2 \cap Z_1}\mathcal{F} \to\Gamma_{Z_2}\mathcal{F} \to \Gamma_{Z_2 \setminus Z_1}\mathcal{F} .
$$

Let $Z$ be a locally closed subset of $X$. We are going to define the functor $i_{Z!}$ such that for  $\F \in {\rm Mod}(k_Z)$, $i_{Z!}\F$ is the unique $k$-sheaf in ${\rm Mod}(k_X)$ inducing $\F$ on $Z$ and zero on $X\setminus Z$. First let $U$ be an open subset of $X$ and let $\F \in {\rm Mod}(k_U)$. Then $i_{U!}\F$ is the sheaf associated to the presheaf $V\mapsto \Gamma(V;i_{U!}\mathcal{F})$ which is $\Gamma(V;\mathcal{F})$ if $V \subseteq U$ and $0$ otherwise. If $S$ is a closed subset of $X$ and $\F \in {\rm Mod}(k_S)$, then $i_{S!}\F=i_{S*}\F$. Now let $Z=U \cap S$ be a locally closed subset of $X$, then one defines $i_{Z!}=i_{U!} \circ i_{S!} \simeq i_{S!} \circ i_{U!}$.
%When $Z$ is a locally closed subset of $X$ one defines the functor
%$i_{Z!}:\mathrm{Mod}(k_Z) \to \mathrm{Mod}(k_X)$ as follows, for each $U$ open in $X$
%$$
%\Gamma(U;i_{Z!}{\mathcal F})=\{ s \in \Gamma(U \cap Z;\mathcal{F})\;;\; {\rm supp}\, s\,\,{\rm is}\,\,{\rm closed}\,\,{\rm in}\,\,X\}
%$$
%i.e.
%$$
%\Gamma(U;i_{Z!}{\mathcal F}) \simeq \lind S \Gamma_S(U \cap Z;\mathcal{F})
%$$
%where $S$ ranges through the family of closed subsets of $X$ contained in $Z$. Of course, when $Z$ is closed $i_{Z!}=i_{Z*}$.
The functor $i_{Z!}$ is exact and has a right adjoint, denoted by $i_Z^!$, when $Z$ is open we have $i_Z^!\simeq i_Z^{-1}$, when $Z$ is closed $i_Z^! \simeq i_Z^{-1}\Gamma_Z$. With these definitions one has
$$
\mathcal{F}_Z \simeq \mathcal{F} \otimes k_Z \simeq i_{Z!}i_Z^{-1}\mathcal{F} \ \ \mathrm{and} \ \ \Gamma_Z\mathcal{F} \simeq \mathcal{H}om(k_Z,\mathcal{F}) \simeq i_{Z*}i_Z^!\mathcal{F}.
$$

Let $X$ be a topological space and  $\Phi $ a family of supports on $X$ (i.e. a collection of closed subsets of $X$ such that: (i)  $\Phi $ is closed under finite unions and (ii) every closed subset of a member of $\Phi $ is in $\Phi$). Recall that for ${\mathcal G} \in \mod(k_X)$, an element $s\in \Gamma (X; {\mathcal G})$ is in $\Gamma _{\Phi }(X;{\mathcal G})$ if and only if its {\it support},
$${\rm supp}\,s=X\setminus \cup \{U\subseteq X:U\,\,{\rm is}\,\,{\rm open}\,\,{\rm in}\,\,X\,\,{\rm and}\,\,s_{|U}=0\},$$
is in $\Phi$, i.e.
$$
\Gamma_\Phi(X;\mathcal{G})=\lind {S \in \Phi}\Gamma_S(X;\mathcal{G}).
$$

The following fact (see \cite{b}, Chaper I, Proposition 6.6) will also be useful later:

\begin{prop} \label{prop section on locally closed}
Let $X$  be a topological spaces, $\Phi$ a  family of supports on $X$, $Z$  a locally closed subset of $X$ and let $i_Z: Z \to X$ be the inclusion.
Let $\mathcal{F}$ be a sheaf in $\mathrm{Mod}(k_Z)$. Then
$$
\Gamma_\Phi(X;i_{Z!}\mathcal{F}) \simeq \Gamma_{\Phi|Z}(Z;\mathcal{F}).
$$
\end{prop}

%\pf
%Let $s \in \Gamma_\Phi(X;i_{W!}\F)$. Then there exists $A \in \Phi$ such that $s \in \Gamma_A(X;i_{W!}\F)$, i.e. $s \in \ker(\Gamma(X;i_{W!}\F) \to \Gamma(X \setminus A;i_{W!}\F))$. By definition it is equivalent to say that $s \in \Gamma(X;\F)$ with $\supp\, s \subset S\subseteq W$, $S$ closed in $X$ and $s_{|X\setminus A}=0$. This means that $s \in \Gamma_{A \cap S}(X;\F)\simeq \Gamma_{A \cap S}(W;\F)$ since $A\cap S$ is closed in $X$. Hence we have
%$$
%\Gamma_A(X;i_{W!}\F) \simeq \lind S\Gamma_{A \cap S}(W;\F),
%$$
%where $S$ ranges through the family of closed subsets of $X$ contained in $W$.
%Remark that all the elements of $\Phi_{|W}$ are of the form $A \cap S$ with $A \in \Phi$ and $S \subset W$ closed in $X$, then passing to the limit of $A \in \Phi$ we obtain the result.
%\qed

\subsection{Sheaves on o-minimal spectral spaces}\label{section shsptop}

Let ${\mathcal M}=(M,<, \ldots )$ be our fixed  arbitrary o-minimal structure.  First observe that in $M$ we have the order topology generated by open definable intervals and in $M^k$ we have the product topology generated by the open boxes. Thus every definable set $X\subseteq M^k$ has the induced topology and we say that a definable subset $Z\subseteq X$ is open (resp. closed) if it is open (resp. closed) with the induced topology. Similarly, we can talk about continuous definable maps $f:X\into Y$ between definable sets. This topology has however a problem: in non-standard o-minimal structures definable sets are usually totally disconnected and never connected or locally compact or compact. So we have to introduce definable analogues of these and other topological notions.

Since we do not want to restrict our work to the affine definable setting, we introduce the notion of definable spaces.  A {\it definable space} is a triple $(X,(X_i,\phi _i)_{i=1}^k)$ where:

 \begin{itemize}
\item[(i)]
$X=\cup \{X_i:i=1,\dots ,k \}$;

\item[(ii)]
each $\phi _i:X_i\into M^{l_i}$ is a bijection such that $\phi _i(X_i)$ is a definable subset of $M^{l_i}$;

\item[(iii)]
for all $j$, $\phi _i(X_i\cap X_j)$ is open in $\phi _i(X_i)$ and  the transition maps $\phi _{ij}:\phi _i(X_i\cap X_j)\into \phi _j(X_i\cap X_j):x\mapsto \phi _j(\phi _i^{-1}(x))$ are definable homeomorphisms.
\end{itemize}

The {\it dimension} of a definable space $X$ is defined as
$$
\dim X=\max\{\dim \phi _i(X_i):i=1,\dots , k\}.
$$
A definable space has a topology such that each $X_i$ is open and the $\phi _i$'s are homeomorphisms: a subset $U$ of $X$ is an open in the basis for this topology if and only if for each $i$, $\phi _i(U\cap X_i)$ is an open definable subset of $\phi _i(X_i)$. We also say that a subset $A$ of $X$ is definable if and only if for each $i$,  $\phi _i(A\cap X_i)$ is a definable subset of $\phi _i(X_i)$.  A map between definable spaces is definable if when it is read through  the charts it is definable. Thus we have the category of definable spaces with definable continuous maps.

We say that  a definable space $X$ is:

 \begin{itemize}
\item
 \textit{definably connected} if it is not the disjoint union of two open and closed definable subsets;

\item
\textit{definably compact} if for every continuous definable map $\sigma :(a,b)\subseteq M\cup \{-\infty, +\infty \}\into X$, the limits $\lim _{t\rightarrow a^+}\sigma (t)$ and  $\lim _{t\rightarrow b^-}\sigma (t)$ exist and belong to $X.$

\item
\textit{definably locally compact} if for every definably compact subset $Z$ with open definable neighborhood $U$ in $X$, there is a definably compact neighborhood of $Z$ in $U$.

\item
\textit{definably normal} if for every disjoint closed definable subsets $Z_1$ and $Z_2$ of $X$ there are disjoint open definable subsets $U_1$ and $U_2$ of $X$ such that $Z_i\subseteq U_i$ for $i=1,2.$

\item
\textit{definable manifold of dimension $n$} if  $\phi _i(X_i)$ is an open definable subset of $M^n$ for every $i=1,\dots ,k.$
\end{itemize}

The {\it o-minimal site} on a definable space $X$ is the category whose objects are open definable subsets of $X$, the morphisms are the inclusions and the admissible covers are  covers by open definable subsets with finite subcoverings. %Thus we have the category of sheaves of abelian groups ${\rm Sh}_{{\rm dtop}}(X)$ on a definable set $X$ equipped with the o-minimal site.

The following results are an easy adaptation of Propositions 6.4.1 and 6.3.3 of \cite{ks2}, replacing ${\mathcal T}_c$ with open definable (indeed we just need the site generated by a family of open subsets closed under finite intersections and whose coverings admit a finite subcover). The first result gives an  easy way to construct o-minimal $k$-sheaves:

\begin{prop}\label{prop psh and sh}
Suppose that $X$ be a definable space. Let ${\mathcal F}$ be a $k$-presheaf on $X$ relative to the o-minimal site on $X$ and assume that:
\begin{enumerate}
\item
${\mathcal F}(\emptyset )=0;$
\item
for any $U$ and  $V$ open definable subsets of $X$ the canonical sequence
$$0\into {\mathcal F}(U\cup V)\into {\mathcal F}(U)\oplus {\mathcal F}(V)\into {\mathcal F}(U\cap V)$$
is exact.
\end{enumerate}
Then ${\mathcal F}$ is a $k$-sheaf on $X$ relative to the o-minimal site on $X$.
\end{prop}

%\pf
%Let $W$ be an open definable subset of $X$ and let  $\{W_1,\dots , W_n\}$ be a finite admissible cover of $W.$ We have to show that the canonical sequence
%$$0\into {\mathcal F}(W)\into \oplus _{1\leq k\leq n}{\mathcal F}(W_k)\into \oplus _{1\leq i<j\leq n}{\mathcal F}(W_{ij})$$
%is exact ($W_{ij}=W_i\cap W_j$). This is enough since the family of such admissible covers is cofinal in the family of admissible covers of $W.$ We show  this claim by induction on $n.$ The result is clear for $n=1$ and it is the hypothesis for $n=2.$ Suppose that the assertion hold for $j\leq n-1$ and set $W'=\bigcup _{1\leq k<n}W_k.$ By the inductive hypothesis the following commutative diagram is exact
%$$
%\xymatrix{ & \,\, & 0 \ar[d] & 0 \ar[d] & \\
%0 \ar[r] & {\mathcal F}(W)  \ar[r] &
%{\mathcal F}(W')\cup {\mathcal F}(W_n)
%\ar[d] \ar[r] & {\mathcal F}(W'\cap W_n) \ar[d]  & \,\, \\
%\,\, & \,\, &
%\bigoplus _{i<n}{\mathcal F}(W_i)\oplus {\mathcal F}(W_n)
%\ar[d] \ar[r] & \bigoplus _{i<n}{\mathcal F}(W_{in})  & \,\, \\
%\,\, & \,\, &
%\bigoplus _{i<j<n}{\mathcal F}(W_{ij}). &  \,\, & \,\, \\
% }
%$$
%Hence the result follows.
%\qed

The second result shows that in this setting the global sections functor commutes with filtrant inductive limits:

\begin{prop}\label{prop sh and lind}
Let $U$ be a open definable subset of $X$ and let $(\F_i)_{i\in I}$ be a filtrant inductive family of sheaves on the o-minimal site associated to $X$. Then
$$
\Gamma(U;\lind i \F_i) \simeq \lind i \Gamma(U;\F_i).
$$
\end{prop}

We  define the {\it o-minimal spectrum} $\tilde{X}$ of a definable space $X$ as in \cite{c}, \cite{cr} and \cite{p}: it is the set of ultrafilters  of definable subsets of $X$. The o-minimal spectrum $\tilde{X}$ of a definable space  $X$ is $T_0$, quasi-compact and a spectral topological space when equipped with the topology generated by the open subsets of the form $\tilde{U}$, where $U$ is an open definable subset of $X$. That is: (i) it has a basis of quasi-compact open subsets, closed under taking finite intersections;
and (ii) each irreducible closed subset is the closure of a unique point.

The {\it dimension} of the o-minimal spectrum $\tilde{X}$ of a definable space $X$ is defined as
$$
\dim \tilde{X}=\dim X.
$$

By a {\it constructible subset of $\tilde{X}$} we mean a subset of the form $\tilde{A}$ where $A$ is a definable subset of $X$.

We also have the {\it o-minimal spectrum $\tilde{f}:\tilde{X}\into \tilde{Y}$ of a continuous definable map} $f:X\into Y$ between definable spaces: given an ultrafilter $\alpha \in \tilde{X}$, $f(\alpha )$ is the ultrafilter in $\tilde{Y}$ determined by the collection $\{A:f^{-1}(A)\in \alpha \}.$

We now recall some results from \cite{ejp} about this tilde functor. Note that these results were stated in \cite{ejp} in the category of  definable sets but are true in the category of definable spaces with exactly the same proofs.

As we saw in \cite{ejp} we have:

\begin{nrmk}\label{tilde on sets and maps}
{\em
The tilde functor is an isomorphism between the boolean algebra of definable subsets of a definable space $X$ and  the boolean algebra of constructible subsets of its o-minimal spectrum $\tilde{X}$ and it commutes with image and inverse image under definable maps.
}
\end{nrmk}

Another useful property is the following result:

\begin{thm}[\cite{ejp}]\label{thm main normal def normal}
Let $X$ be a definable space. Then  the following hold:
\begin{itemize}
\item[(i)]
$X$ is definably  connected if and only if its  o-minimal spectrum $\tilde{X}$ is connected.
\item[(ii)]
$X$ is definably normal if and only if its  o-minimal spectrum $\tilde{X}$ is normal.
\end{itemize}
\end{thm}

Also we have the following shrinking lemma:

\begin{prop}[\cite{ejp}, The shrinking lemma]\label{prop shrinking lemma}
Suppose that $X$ is a definably normal definable space (resp. a normal o-minimal spectrum of a definable space). If $\{U_i:i=1,\dots ,n\}$ is a covering of $X$ by open definable subsets (resp. open subsets) of
$X$, then there are definable (resp. constructible) open subsets $V_i$ and definable (resp.
constructible) closed subsets $K_i$ of $X$ ($1\leq i\leq n$) with $V_i\subseteq K_i\subseteq U_i$ and $X=\cup \{V_i:i=1,\dots, n\}$.
\end{prop}

Since the o-minimal spectrum of a definable space is quasi-compact, as in the proof of Propositions \ref{prop psh and sh} and \ref{prop sh and lind}, we have:

\begin{nrmk}\label{nrmk psh and sh spec}
{\em
Suppose that $X$ is an object in the category of  o-minimal spectra of definable spaces. Let ${\mathcal F}$ be a $k$-presheaf on $X$ and assume that:
\begin{enumerate}
\item
${\mathcal F}(\emptyset )=0;$
\item
for any $U$ and  $V$ open constructible subsets of $X$ the canonical sequence
$$0\into {\mathcal F}(U\cup V)\into {\mathcal F}(U)\oplus {\mathcal F}(V)\into {\mathcal F}(U\cap V)$$
is exact.
\end{enumerate}
Then ${\mathcal F}$ is a $k$-sheaf on $X$. Moreover sections on open constructible subsets commute with filtrant $\Lind$.
}
\end{nrmk}

We have a morphism of sites naturally induced by the above tilde functor from the category of definable spaces with continuous definable maps into the category of o-minimal spectral spaces with the o-minimal spectra of continuous definable maps. This morphism of sites induces the following isomorphism:

 \begin{thm}[\cite{ejp}]\label{thm main iso on sheaves}
Let $X$ be a definable space. Then there is an isomorphism between the category of $k$-sheaves  on $X$ relative to o-minimal site on $X$ and the category of $k$-sheaves  on the o-minimal spectrum $\tilde{X}$ of $X$  relative to the spectral topology on  $\tilde{X}.$
\end{thm}

 The isomorphism of Theorem \ref{thm main iso on sheaves} allowed the development of o-minimal sheaf cohomology without supports in \cite{ejp} by defining concepts and also sometimes proving results via this tilde isomorphism.  In this paper we will continue to use this technique but allowing now the presence of supports.

%Hence we  define  sheaf cohomology in the category of pairs of definable spaces by setting $$H^*(X;{\mathcal F}):=H^*(\tilde{X};\tilde{{\mathcal F}})$$ where $X$ is a definable set, ${\mathcal F}$ is a sheaf in ${\rm Sh}_{{\rm dtop}}(X)$. (In this context we  have  \cite{ejp} Proposition 3.2 giving an isomorphism of categories $${\rm Sh}_{{\rm dtop}}(X)\into {\rm Sh}(\tilde{X}):{\mathcal F}\mapsto \tilde{{\mathcal F}},$$ where ${\rm Sh}(\tilde{X})$ is the category of sheaves of abelian groups on the topological space $\tilde{X}$).

We will now define the notion of family of supports on a definable set. Our treatment of this will follow the corresponding theory in semi-algebraic geometry in \cite{D3} (Chapter II, Sections 1 - 5) and in topological spaces in \cite{b} (Chapter I, Section 6 and Chapter II, Sections 9 - 13). Note also that since, as we saw in \cite{ejp}, the role of paracompactness in sheaf theory on topological spaces has to be replaced by normality in sheaf theory on o-minimal spectral spaces, we will continue to do this here.

\begin{defn}\label{defn def supports}
{\em
Let $X$ be a definable space. A {\it family of definable supports} is a  family of closed definable subsets of $X$ such that:

\begin{enumerate}
\item
every closed definable subset of a member of $\Phi$ is in $\Phi $;

\item
$\Phi $ is closed under finite unions.

\vspace{.15in}
\noindent
$\Phi  $ is said to be a {\it family of  definably normal supports} if in addition:

\item
each element of $\Phi $ is definably normal;

\item
for each element $S$ of $\Phi  $, if $U$ is an open definable neighborhood of $S$ in $X$, then there exists a (closed) definable neighborhood of $S$ in $U$ which is in $\Phi$.
\end{enumerate}
}
\end{defn}

\begin{expl}\label{expl supp c}
{\em
Let $X$ be a definable space and let $c$ be the collection of all definably compact definable subsets of $X.$ Then $c$ is a family of definable supports on $X.$ Moreover, if $X$ is definably normal and definably locally compact, then $c$ will be a family of  definably normal supports on $X.$
}
\end{expl}

If $Y$ is a definable subset of the definable space $X$ and $\Phi  $ a family of definable supports on $X$, then we have families of definable  supports
$$\Phi  \cap Y=\{A\cap Y:A\in \Phi \}$$
and
$$\Phi _{|Y}=\{A\in \Phi : A\subseteq Y\}$$
on $Y$.

If $f:X\into Z$ is a continuous definable map between definable spaces and $\Phi $ is a family of definable supports on $Z$,
then we have a family of definable supports
$$f^{-1}\Phi  =\{A\subseteq X: A\,\,{\rm is} \,\, {\rm closed,}\,\,{\rm definable}\,\,{\rm and}\,\,\exists B\in \Phi \,\,(A\subseteq f^{-1}(B)\}$$
on $X$.

%\begin{defn}\label{defn spec supports}
%{\em
%Let $X$ be an object in the category of o-minimal spectra of  definable space. A family of supports is a  family of closed  subsets of $X$ such that:

%\begin{enumerate}
%\item
%every closed  subset of a member of $\Phi$ is in $\Phi $;

%\item
%$\Phi $ is closed under finite unions.

%\vspace{.15in}
%\noindent
%$\Phi  $ is said to be a {\it family of  normal supports} if in addition:

%\item
%each element of $\Phi $ is  normal;

%\item
%for each element $S$ of $\Phi  $, if $U$ is an open  neighborhood of $S$ in $X$, then there exists a (closed) constructible  neighborhood of $S$ in $U$ which is in $\Phi$.
%\end{enumerate}
%}
%\end{defn}

%A
\begin{nrmk}\label{nrmk def supp and cons supp}
{\em
Note that a family of definable supports $\Phi $ on a definable space $X$ determines a family of supports
$$\tilde{\Phi }=\{A\subseteq \tilde{X}: A\,\,{\rm is} \,\,{\rm closed}\,\,{\rm and}\,\,\exists B\in \Phi \,\,(A\subseteq \tilde{B})\}$$
on the topological space $\tilde{X}$.  By Remark \ref{tilde on sets and maps}  it follows that
$$\tilde{\Phi  \cap Y}=\tilde{\Phi }\cap \tilde{Y},\,\, \tilde{\Phi  _{|Y}}=\tilde{\Phi }_{|\tilde{Y}}\,\,{\rm and}\,\,\tilde{f^{-1}\Phi }=\tilde{f}^{-1}\tilde{\Phi }.$$
We will say that the family of supports on $\tilde{X}$ is {\it constructible}  if it is obtained by applying tilde to some family of definable supports on $X$.

By theorem \ref{thm main normal def normal} it follows that  $\Phi $ is definably normal if and only if $\tilde{\Phi}$ is normal. Here, we say that $\Psi $ is  a {\it family of  normal supports} on the spectral topological space $\tilde{X}$  if is a family of supports and:
\begin{enumerate}
\item
each element of $\Psi $ is  normal;

\item
for each element $S$ of $\Phi  $, if $U$ is an open  neighborhood of $S$ in $\tilde{X}$, then there exists a (closed) constructible  neighborhood of $S$ in $U$ which is in $\Phi$.
\end{enumerate}
}
\end{nrmk}

\begin{defn}\label{defn  sheaf cohomo Phi}
{\em
Let $X$ be a definable space, $\Phi $ a family of definable supports in $X$  and ${\mathcal F}$ a $k$-sheaf on $X$ relative to the o-minimal site on $X$. We define the {\it o-minimal sheaf cohomology groups with definable supports in $\Phi$ } via the tilde isomorphism of Theorem \ref{thm main iso on sheaves} by
$$H^*_{\Phi}(X;{\mathcal F})=H^*_{\tilde{\Phi}}(\tilde{X};\tilde{{\mathcal F}}),$$
where $\tilde{\mathcal{F}}$ is the image of $\mathcal{F}$ via the isomorphism between the category of $k$-sheaves  on $X$ relative to o-minimal site on $X$ and the category of $k$-sheaves  on the o-minimal spectrum $\tilde{X}$ of $X$.\\

If $f:X\into Y$ is a continuous definable map, we define the induced homomorphism
$$f^*:H^*_{\Phi }(Y;{\mathcal F})\into H^*_{f^{-1}\Phi }(X;f^{-1}{\mathcal F})$$
in cohomology to be the
same as the induced homomorphism
$$\tilde{f}^*:H^*_{\tilde{\Phi }}(\tilde{Y};{\tilde {\mathcal F}})\into
H^*_{\tilde{f}^{-1}\tilde{\Phi }}(\tilde{X};\tilde{f}^{-1}{\tilde {\mathcal F}})$$
in cohomology of the continuous map
$\tilde{f}:\tilde{X}\into \tilde{Y}$ of topological spaces.
}
\end{defn}

The proof of the o-minimal Vietoris-Begle theorem with supports below is similar to its analogue  without supports (\cite{ejp} Theorem 4.3) using classical arguments:

\begin{thm}[Vietoris-Begle theorem]\label{thm vietoris-b}
Let $f:X\into Y$ be a surjective  morphism  in the category of o-minimal spectra of definable spaces  that maps constructible closed subsets of $X$ onto  closed subsets of $Y$. Let ${\mathcal F}\in \mod (k_Y)$, $\Phi $ a constructible family of supports on $Y$  and suppose that $Y$ is a subspace of a normal space  in the category of o-minimal spectra of definable spaces. Assume that $f^{-1}({\beta  })$ is connected and $H^q(f^{-1}({\beta  });
f^{-1}{\mathcal F}_{|f^{-1}({\beta  })})=0$ for $q>0$ and all ${\beta  }\in Y$. Then the induced map
$$f^*:H^*_{\Phi }(Y;{\mathcal F})\into H^*_{f^{-1}\Phi }(X;f^{-1}{\mathcal F})$$
is an isomorphism.
\end{thm}

%\begin{defn}\label{defn rel sheaf cohomo c}
%{\em
%If $X$ is a definable set, $A$ is a closed definable subset of $X$  and ${\mathcal F}$ a sheaf in ${\rm Sh}_{{\rm dtop}}(X)$, we define as above the {\it relative o-minimal sheaf cohomology groups with definably compact supports} $$H^*_{c }(X,A;{\mathcal F})$$ by replacing ${\mathcal F}$ by ${\mathcal F}_{X\setminus A}$. Similarly, if $f:(X,A)\into (Y,B)$ is a continuous definable map of closed pairs of definable sets (i.e., $A\subseteq X$ and $B\subseteq Y$ are closed definable subsets and $f:X\into Y$ is a continuous definable map such that $f(A)\subseteq B$) and  ${\mathcal F}$ a sheaf in ${\rm Sh}_{{\rm dtop}}(Y)$, then the induced homomorphisms $$f^*:H^*_{c }(Y,B;{\mathcal F})\into H^*_{f^{-1}c }(X,A;f^*{\mathcal F})$$ in cohomology are defined as above by replacing ${\mathcal F}$ by ${\mathcal F}_{Y\setminus B}$.
%}
%\end{defn}

We have in this context  the Eilenberg-Steenrod axioms with definable supports adapted to the o-minimal site. Indeed, once we pass to the category of o-minimal spectra of definable spaces the proofs of the exactness and excision axioms are  purely algebraic. See \cite{b}. The dimension axiom is also immeadiate. On the other hand, from the  Vietoris-Begle theorem (Theorem \ref{thm vietoris-b}) we obtain:

\begin{thm}[Homotopy axiom]\label{thm homotopy ax scs}
Suppose that $X$ is a  definable space  and ${\mathcal F}$ is a $k$-sheaf on $X$ relative to the o-minimal site on $X$. Let $[a,b]\subseteq M$ be a closed interval. Assume that $\M$ has definable Skolem functions, $X$ is definably normal and the projection $ \pi:X\times [a,b]\into X$ maps closed definable subsets of $X\times [a,b]$ onto closed definable subsets of $X$. If for $d\in [a,b]$,
$$i_d:X\into X\times [a,b]$$
is the continuous definable map given by $i_d(x)=(x,d)$ for all $x\in X$, then
\[ i_a^*=i_b^*:H^n_{\Phi  \times [a,b]}(X\times [a,b];\pi ^{-1}{\mathcal F})\into
H^n_{\Phi  }(X;{\mathcal F}) \]
for all $n\in \NN$.
\end{thm}

\pf
The homotopy axiom will follow once we show that the projection map
$\pi :X\times [a,b]\into X$ induces an isomorphism
$$\pi ^*:H^n_{\Phi  }(X;{\mathcal F})\into H^n_{\Phi  \times [a,b]}(X\times [a,b];\pi ^{-1}{\mathcal F})$$
since by functoriality we obtain
$$i_a^*=i_b^*=(\pi ^*)^{-1}:H^n_{\Phi  \times [a,b]}(X\times [a,b];\pi ^{-1}{\mathcal F})
\into H^n_{\Phi  }(X;{\mathcal F})$$
for all $n\in \NN$.  Equivalently we
need to show that
$$\tilde{\pi }^*:H^n_{\tilde{\Phi  }}(\tilde{X};{\tilde {\mathcal F}})\into H^n_{\tilde{\Phi  \times [a,b]}}(\tilde{X\times [a,b]};\tilde{\pi }^{-1}{\tilde {\mathcal F}})$$
is an isomorphism. For
this we need to verify the hypothesis of the Vietoris-Begle theorem (Theorem \ref{thm vietoris-b}), but this was done in the proof of the homotopy axiom for o-minimal sheaf cohomology without supports
(\cite{ejp} Theorem 4.4).
\qed

\begin{nrmk}\label{nrmk exact triple mv}
{\em
In this context we also have the exactness for triples of closed definable subsets  and the Mayer-Vietoris theorem for $\Phi $-excisive pairs of definable sets. See \cite{b}.
}
\end{nrmk}

\end{section}

\begin{section}{$\Phi$-soft sheaves}\label{section appendix}

The results we present below are in the category of o-minimal spectra of definable spaces but by the
isomorphism of Theorem \ref{thm main iso on sheaves}  they have a suitable, but more
restrictive, analogue in the category of definable spaces. In fact these results are the analogue of classical results on paracompactifying families of supports on topological spaces (\cite{b}) adapted to normal and constructible families of supports on spectral spaces.

\begin{subsection}{Normal and constructible supports}\label{subsection normal and constructible supports}

We start the subsection with the following useful result:

%We have the chain of isomorphisms
%\begin{eqnarray*}
%\Gamma_\Phi(X;i_{W!}\mathcal{F}) & \simeq & \Gamma_\Phi(X;\lind A i_{W*}\Gamma_A\mathcal{F}) \\
%& \simeq & \lind A \Gamma_\Phi(X; i_{W*}\Gamma_A\mathcal{F}) \\
%& \simeq & \lind A \Gamma_{\Phi\cap A}(W;\mathcal{F}) \\
%& \simeq & \Gamma_{\Phi|W}(W;\mathcal{F}) ,
%\end{eqnarray*}
%where $A$ ranges through the family of closed subsets of $X$ contained in $W$. In the second isomorphism we used Proposition \ref{prop cohomo limits soft}, in the second one the isomorphism $i_{W*}\Gamma_A \simeq \Gamma_Ai_{W*}$, and the last one follows since every element of $\Phi_{|W}$ is the intersection of a closed subset of $W$ contained in $A$ and an element of $\Phi$.
%\qed

\begin{prop}\label{prop cohomo limits soft}
Assume that $X$ is an object in the category of o-minimal spectra of definable spaces  and let $({\mathcal F}_i)_{i\in I}$ be a filtrant inductive family of sheaves in $\mathrm{Mod}(k_X)$ and $\Phi $ a  constructible family of supports on $X$. Then
$$ \Gamma_{\Phi }(X;\varinjlim_{i\in I}{\mathcal F}_i)=\varinjlim_{i\in I}\Gamma _{\Phi }(X;{\mathcal F}_i).$$
\end{prop}

\pf
First observe that by definition for any $\G\in\mod(k_X)$ we have
$$\Gamma _{\Phi }(X; {\mathcal G})=\varinjlim_{S\in \Phi}\Gamma _{S}(X;{\mathcal G}).$$
Thus it is enough to show that for each $S\in \Phi $ constructible we have
$$\Gamma_{S}(X;\varinjlim_{i\in I}{\mathcal F}_i)=\varinjlim_{i\in I}\Gamma _{S}(X;{\mathcal F}_i).$$
For this consider the following commutative diagram where the vertical arrows are the canonical maps
{\tiny
\[\xymatrix{
0\ar[r] & \lind {i\in I}\Gamma _S(X;{\mathcal F}_i) \ar[d]\ar[r] &  \lind {i\in I}\Gamma (X;{\mathcal F}_i) \ar[d]\ar[r] &  \lind {i\in I}\Gamma (X\setminus S; {\mathcal F}_i)\ar[d] \\
0 \ar[r] & \Gamma _S(X;\lind {i\in I}{\mathcal F}_i) \ar[r] & \Gamma (X;\lind {i\in I}{\mathcal F}_i)  \ar[r] &\Gamma (X\setminus S;\lind {i\in I}{\mathcal F}_i)
.}
\]}
The rows are exact  by definition of $\Gamma_S$ (i.e. $\Gamma_S(X;\G)=\ker(\Gamma(X;\G)\to\Gamma(X\setminus S;\G)$ for any $\G\in\mod(k_X)$) and by the exactness of filtrant $\Lind$.
Since  $X$ and $X\setminus S$ are  open constructible subsets of $X$ and sections on open constructible subsets commute with filtrant $\Lind$ (Remark \ref{nrmk psh and sh spec}),  it follows that the two vertical arrows on the right are isomorphisms. Hence  the first vertical arrow is also an isomorphism as required.
\qed

%\begin{clm}\label{clm cons  and limits}
%If $U$ is an open constructible subset of $X$, then
%$$\Gamma (U,\lind {i\in I}{\mathcal F}_i)=\lind {i\in I}\Gamma (U,{\mathcal F}_i).$$
%\end{clm}

%\pf
% Denote by $\indl {i\in I}{\mathcal F}_i$ the presheaf $V\mapsto \lind {i\in I}\Gamma (V,{\mathcal F}_i)$ on $X.$ Let $U$ be an open constructible subset of $X$ and let ${\mathcal S}\in \cov (U)$ be a finite admissible cover of $U$. Since $\lind {i}$ commutes with finite projective limits we have an isomorphism
% $$(\indl {i}{\mathcal F}_i)({\mathcal S})\simeq \lind {i}{\mathcal F}_i({\mathcal S}).$$
%Since the family of  finite admissible covers ${\mathcal S}\in \cov (U)$ of $U$  is cofinal in the family of admissible covers $\cov (U)$ of $U$, because $U$ is quasi-compact,  and ${\mathcal F}_i(U)\simeq {\mathcal F}_i({\mathcal S})$ since ${\mathcal F}_i\in \mod (k_X)$  for each $i,$ we have
%$$\indl {i}{\mathcal F}_i\simeq (\indl {i}{\mathcal F}_i)^+.$$
%Applying again the  sheafification functor $(\bullet )^+$ we get
%$$\indl {i}{\mathcal F}_i\simeq (\indl {i}{\mathcal F}_i)^+\simeq (\indl {i}{\mathcal F}_i)^{++}\simeq \lind {i}{\mathcal F}_i$$
%and the result follows by applying the functor $\Gamma (U, \bullet ).$

% Since the spectral topology on $X$ (which is the topology generated by the open constructible subsets) is noetherian, i.e., the open constructible subsets are quasi-compact, the result follows at once from \cite{a}, Chapter II, 5.3.
%\qed

The following lemma is fundamental in this Subsection:

\begin{lem}\label{clm extending sections no supports}
Assume that  $Z$ is a subspace of a normal space $X$ in the category of o-minimal spectra of definable spaces, ${\mathcal G}$ is a sheaf in $\mathrm{Mod}(k_Z)$ and $Y$ is a quasi-compact subset of $Z$. Then the canonical morphism
$$\lind {Y\subseteq U}\Gamma(U\cap Z;{\mathcal G})
\into \Gamma(Y;{\mathcal G}_{|Y})\,\,\,$$
where $U$ ranges through the family of open constructible subsets of $X$, is an isomorphism.
\end{lem}

\pf
%Let $X$ be the normal o-minimal spectrum of a definable space of which $Z$ is a subspace.
Since $Y$ is quasi-compact, the family of open neighborhoods of $Y$ in  $Z$ of the form $V\cap Z$ where $V$ is an open constructible subset of $X$ is a fundamental system  of neighborhoods of $Y$ in $Z$. Hence, the morphism of the lemma is certainly injective.

To prove that it is surjective, consider a section $s\in \Gamma (Y;\G _{|Y})$.
There is a covering $\{U_j:j\in J\}$ of $Y$ by open constructible subsets of
$X$ and sections $s_j\in \Gamma (U_j\cap Z; \G_{|U_j\cap Y})$, $j\in J$, such
that $s_{j|U_j\cap Y}=s_{|U_j\cap Y}$. Since $Y$ is quasi-compact, we can
assume that $J$ is finite, and so $\cup \{U_j:j\in J\}$ is an open
constructible subset of $X$.
Since $X$ is normal, by the shrinking lemma
(Proposition \ref{prop shrinking lemma}), there are open constructible subsets
$\{V_j:j\in J\}$ of this union such that $\overline{V}_j\subseteq U_j$ for
every $j\in J$ and $Y\subseteq \cup \{V_j:j\in J\}.$ For $x \in Z$ set $J(x)=\{j\in J:x\in \overline{V}_j\}$. Each $x$ has a constructible neighborhood $W_x$ with $J(y)
\subseteq J(x)$ for each $y \in W_x.$ This is defined by
$$
W_x = (\bigcap_{x
\in V_l}V_l \cap \bigcap_{j\in J(x)}U_j) \setminus \bigcup_{k \notin J
(x)}\overline{V}_k.
$$
Observe that for all $i,j \in J(x)$ we have that $W_x$ is an open subset of both $U_i$ and  $U_j$. Hence, for every $i,j \in J(x)$ we have $s_{i|W_x\cap Y}= s_{|W_x\cap Y}=s_{j|W_x\cap Y}$. So, for  $y \in W_x \cap Y$, we have  $(s_i)_y=(s_j)_y$ for any $i,j \in J
(x)$.
This implies that the set
$$W=\{x\in (\bigcup_{j\in J}V_j)\cap Z:\textrm
{$(s_i)_x=(s_j)_x$ for any $i,j\in J(x)$}\}$$
contains $Y$ (clearly $Y\subseteq \bigcup _{x\in Z} W_x\cap Y\subseteq (\bigcup_{j\in J}V_j)\cap Z$). On the other hand,  the condition $(s_i)_z=(s_j)_z$ for any $i,j \in J
(x)$ and the fact that $J(x)$ is finite implies that $z$ has an open  neighborhood in $Z$  on which $s_i=s_j$ for any $i,j \in J(x)$. Thus $W$ is an open neighborhood of $Y$ in $Z$.
Since $Y$ is quasi-compact we may assume that  $W$ is of the form $U\cap Z$ for
some open constructible subset $U$ of $X$. Since $s_{i_{|W \cap V_i \cap V_j}}
=s_{j_{|W \cap V_i \cap V_j}}$
there exists $t\in \Gamma (W;\G )$ such that $t_
{|W \cap V_j}=s_{j|W \cap V_j}$. This proves that the morphism is surjective.
\qed

%As we pointed out in \cite{ejp} page 173, the proof of  Lemma \ref{clm extending sections no supports} is a consequence of the shrinking lemma (Proposition \ref{prop shrinking lemma}) and holds in any normal spectral space. For details see its analogue in \cite{D2} Lemma 2.2.

A general form of Lemma \ref{clm extending sections no supports} is:

\begin{lem}\label{lem extending sections}
Assume that $X$ is an object in the category of o-minimal spectra of definable spaces, $Z$ is a subspace of $X$, ${\mathcal G}$ is a sheaf in ${\rm Mod}(k_Z)$, $\Phi $ is a normal and constructible family of supports on $X$ and $Y$ is a subset of $Z$ such that $D\cap Y$ is a quasi-compact  subset for every $D\in \Phi $.
%Then for every $s\in \Gamma _{\Phi \cap Y}(Y,{\mathcal F}_{|Y})$ there exists an open neighborhood $W$ of $Y$ in $Z$ and $t\in \Gamma _{\Phi \cap W}(W,{\mathcal F}_{|W})$ such that $t_{|Y}=s$.
Then the canonical morphism
$$\lind {Y\subseteq U}\Gamma_{\Phi \cap U \cap Z}(U\cap Z;{\mathcal G})
\into \Gamma_{\Phi \cap Y}(Y;{\mathcal G}_{|Y})\,\,\,$$
where $U$ ranges through the family of open constructible subsets of $X$, is an isomorphism.
\end{lem}

\pf
Let us prove injectivity. Let $s \in \Gamma_{D\cap U \cap Z}(U\cap Z;{\mathcal G})$, with $D \in \Phi$ and $U \supset Y$ open constructible subset of $X$ and such that $s_{|D \cap Y}=0$. Since $\Phi $ is a normal and constructible family of supports on $X$, there is a constructible and normal $E\in \Phi $  which is a closed neighborhood of  $D$ in $X$. Thus $D\cap Z$ is a subspace of a normal space $E$ in the category of o-minimal spectra of definable spaces and $D\cap Y$ is a quasi-compact subset of $D\cap Z$. By Lemma \ref{clm extending sections no supports} applied to $E$, $D\cap Z$ and $D \cap Y$, there exists an open (in $E$) constructible neighborhood $V'$ of $D\cap Y$ such that $s_{|V'\cap D\cap Z}=0$. Of course we may assume that $V'=V\cap E$ for some open constructible subset $V$ of $X$. So there exists an open (in $X$) constructible neighborhood $V$ of $D\cap Y$ such that $s_{|V\cap D\cap Z}=0$.
Also, by replacing $V$ with its intersection with $U$ if necessary we may assume that  $V \subseteq U.$
 %of $D \cap Y$ such that $s_{|V\cap Z}=0$.
 Set $W=V \cup(U \setminus D)$. Then  $W$ is open constructible in $X$, $Y\subseteq W\subseteq U$ and $s_{|W\cap Z}=0$.

Let us prove that the morphism is surjective. Let $s\in \Gamma _{\Phi \cap Y}(Y;
{\mathcal G}_{|Y})$ and consider normal constructible sets $C$, $D$  and
$E$ in $\Phi $ such that $D$ is a closed neighborhood of $C$ in $X$,  $E$ is a
closed neighborhood of $D$ in $X$ and the support of $s$ is contained in $C\cap Y$.
We shall find $\widetilde{t} \in \Gamma_D(U \cap Z;{\mathcal G})$ such that
$\widetilde{t}_{|Y}=s$.
%Replacing $X$ by $E \setminus \partial E$ we may assume that $X$ is normal.
After applying Lemma \ref{clm extending sections no
supports} above to $E$, $D\cap Z$ and $D \cap Y$
%and using quasi-compactness of $D\cap Y$,
we see that there exists an open in $E\setminus \partial E$ (and hence in $X$) constructible neighborhood $V$ of $D\cap Y$ and a section
$t\in \Gamma (V\cap D\cap Z;{\mathcal G})$ such that $t_{|D\cap Y}=s_{|D\cap Y}$.
Since $t_{|\partial D\cap Y}=0$, then each point $x$ of $\partial D \cap Y$ has
an open constructible neighborhood $W_x \subset V$ such that $t_{|W_x\cap D\cap Z}=0.$ Using
quasi-compactness of $\partial D\cap Y$ (it is closed on the quasi-compact set $D\cap Y$), there exists a finite number of points
$x_1,\ldots,x_n$ such that $\partial D \cap Y \subset \bigcup_{i=1}^nW_{x_i}:
=W$. We have $t_{|W\cap D \cap Z}=0$ and $W$ is open constructible. Let $U_1=(V\cap
(D\setminus \partial D))\cup W$. Then $U_1$ is open constructible and $D\cap
Y\subseteq U_1\subseteq V$. Define $t'«\in \Gamma (U_1\cap Z; {\mathcal G})$ by: $t'_{|V\cap (D\setminus \partial D)\cap Z}=t_{|V\cap (D\setminus \partial D)\cap Z}$ and $t'_{|W\cap Z}=0$. This is well defined since $t_{|W\cap D \cap Z}=0$ and $(V\cap (D\setminus \partial D)\cap Z)\cap (W\cap Z)\subseteq W\cap D \cap Z.$ Observe also that $t'_{|U_1\cap D\cap Z}=t.$ Let $U_2=X\setminus D$. Then $U=U_1\cup U_2$ is
open constructible, $Y\subseteq U$, $U_1\cap U_2\subseteq W$ and we can define
$\widetilde{t} \in \Gamma(U\cap Z;\mathcal{G})$ in the following way:
$\widetilde{t}_{|U_1\cap Z}=t'_{|U_1\cap Z}$, $\widetilde{t}_{|U_2\cap Z}=0$. It
is well defined since $t'_{|W\cap Z}=0$ and $U_1 \cap U_2 \subset W$. Moreover
$\mathrm{supp}\,\widetilde{t} \subseteq D$ and $\widetilde{t}_{|Y}=s$ as
required.
\qed

Recall that a  sheaf ${\mathcal F}$ on a topological space $X$ with a family of supports $\Phi $ is {\it $\Phi $-soft} if and only if  the restriction $\Gamma (X;{\mathcal F})\into \Gamma (S;{\mathcal F}_{|S})$ is surjective for every   $S\in \Phi .$
%the restriction $\Gamma _{\Phi }(X;{\mathcal F})\into \Gamma _{\Phi \cap Y}(Y;{\mathcal F}_{|Y})$ is surjective.  Since for all $Y\in \Phi$ we have $\Phi \cap Y=\Phi _{|Y}$, the sheaf ${\mathcal F}$ on the  topological space $X$ is $\Phi $-soft if and only if for every  subset $Y$ of $X$ which is in $\Phi $  the restriction $\Gamma _{\Phi }(X;{\mathcal F})\into \Gamma _{\Phi _{|Y}}(Y;{\mathcal F}_{|Y})$ is surjective.
If $\Phi $ consists of all closed subsets of $X$, then ${\mathcal F}$ is simply called {\it soft}.

\begin{prop}\label{prop phi soft}
Let $X$ be a topological space and  ${\mathcal F}$ is a sheaf in ${\rm Mod}(k_X)$. If  $\Phi $  is a family of supports on $X$ such that  every $C\in \Phi$ has a neighborhood $D$ in $X$ with $D\in \Phi$. Then the following are equivalent:
\begin{enumerate}
\item
 ${\mathcal F}$ is $\Phi $-soft;
 \item
 ${\mathcal F}_{|S}$ is soft for every  $S\in \Phi $;
  \item
 $\Gamma _{\Phi }(X;{\mathcal F})\into \Gamma _{\Phi _{|S}}(S; {\mathcal F}_{|S})$ is surjective  for every closed subset $S$ of $X$;

 \medskip
 \noindent
If  in addition $X$ is an object  in the category of o-minimal spectra of definable spaces and $\Phi $ is a  constructible family of supports on $X$, then the above are also equivalent to:

\medskip
 \item
 ${\mathcal F}_{|Z}$ is soft for every constructible subset $Z$ of $X$ which is in $\Phi $;

 \medskip
 \noindent
If  moreover  $\Phi $ is a  normal and constructible family of supports on $X$, then the above are also equivalent to:

\medskip

  \item
 $\Gamma (X;{\mathcal F})\into \Gamma (Z; {\mathcal F}_{|Z})$ is surjective  for every constructible subset $Z$ of $X$ which is in $\Phi $;
%  \item
 %$\Gamma _{\Phi }(X;{\mathcal F})\into \Gamma _{\Phi \cap Z}(Z;{\mathcal F}_{|Z})$ is surjective  for every constructible subset $Z$ of $X$ which is in $\Phi $.
  %\item
 %$\Gamma _{\Phi }(X;{\mathcal F})\into \Gamma _{\Phi \cap S}(S;{\mathcal F}_{|S})$ is surjective  for every closed constructible subset $S$ of $X$.
\end{enumerate}
\end{prop}

\pf
The equivalence of  (1), (2) and (3) is shown in \cite{b} Chapter II, 9.3. (Our hypothesis is sufficient in the proof given there). The equivalence of (2) and (4) is obvious since every $S\in \Phi $ is contained in some constructible subset  of $X$ which is in $\Phi $.

Clearly (1) implies (5). Assume (5) and let $S\in \Phi $ and $s\in \Gamma (S; {\mathcal F}_{|S})$. Since $\Phi $ is normal and constructible, there is a normal closed  and constructible neighborhood $D$ of $S$ which is in $\Phi $. By Lemma \ref{clm extending sections no supports}, $s$ can be extended to a section $t\in \Gamma (W;{\mathcal F})$ of ${\mathcal F}$ over a neighborhood $W$ of $S$ in $D$. Applying  the shrinking lemma
%(\cite{ejp} Proposition 2.17),
we find a closed constructible neighborhood $Z$ of $S$ in $W$. Since $D\in \Phi $ we have $Z\in \Phi $. So $t_{|Z}\in \Gamma (Z;{\mathcal F}_{|Z})$, $(t_{|Z})_{|S}=s$ and $t_{|Z}$ can be extended to $X$ by (5). Hence, (5) implies (1).

\qed

\begin{cor} \label{cor lim soft}
Assume that $X$ is an object  in the category of o-minimal spectra of definable spaces and $\Phi $ is a normal and constructible family of supports on $X$. Then filtrant inductive limits of $\Phi$-soft sheaves in ${\rm Mod}(k_X)$  are $\Phi$-soft.
\end{cor}

\pf
It follows combining Propositions \ref{prop cohomo limits soft}  and  \ref{prop phi soft} (5) and the exactness of filtrant inductive limits.
\qed

The following topological result will also be useful below:

\begin{prop}\label{prop soft loc and phi-dim}
Let $X$ be a topological space and  $\Phi $  is a family of supports on $X$ such that  every $C\in \Phi$ has a neighborhood $D$ in $X$ with $D\in \Phi$.
%Assume that $X$ is an object  in the category of o-minimal spectra of definable spaces and  $\Phi $ is a normal and  constructible family of supports on $X$.
Let $W$ be a locally closed subset of $X$. The following hold:

\begin{itemize}
\item[(i)]
if ${\mathcal F} \in {\rm Mod}(k_X)$ is $\Phi$-soft, then $\F_{|W}$ is $\Phi_{|W}$-soft.

\item[(ii)]
${\mathcal G}$ in ${\rm Mod}(k_W)$ is $\Phi _{|W}$-soft if and only if  $i_{W!}{\mathcal G}$ is $\Phi $-soft.

\item[(iii)]
if ${\mathcal F}\in {\rm Mod}(k_X)$ is $\Phi $-soft, then ${\mathcal F}_W$ is $\Phi $-soft.

\end{itemize}
\end{prop}

\pf
(i) If $W$ is open it is obvious. If $W$ is closed it follows from Proposition \ref{prop phi soft} (3). Combining these two cases (i) follows.

(ii) The ``if '' part follows from Proposition \ref{prop phi soft} (2). For the ``only if'' part note that by Proposition  \ref{prop section on locally closed}  (applied to $X$ and $S$ respectively and using the fact that $(\Phi _{|S})_{|W}=\Phi _{|W\cap S}$) we have  $\Gamma_{\Phi}(X;i_{W!}{\mathcal G}) \simeq \Gamma_{\Phi|W}(W;{\mathcal G})$ and $\Gamma _{\Phi |S}(S;i_{W!}{\mathcal G})=\Gamma_{\Phi|W\cap S}(S \cap W;{\mathcal G})$ for any closed subset $S$ of $W$. Then apply Proposition \ref{prop phi soft} (3).

(iii) The result follows from (i) and (ii), since ${\mathcal F}_W=i_{W!}{\mathcal F}_{|W}$.
\qed

%(i)  First observe that we have the following general natural isomorphisms  $\Gamma _{\Phi }(X,{\mathcal G}^X)\simeq \Gamma _{\Phi |W}(W,{\mathcal G})$(\cite{b} Chapter I, 6.6) and similarly, for $D\subset X$ closed, $\Gamma _{\Phi |D}(D,({\mathcal G}^X)_{|D})\simeq \Gamma _{\Phi |W\cap D}(W\cap D,{\mathcal G}_{|W\cap D}).$

%Suppose that ${\mathcal G}$ is $\Phi _{|W}$-soft. If $W=O\cap C$ with $C$ (resp. $O$) a constructible  closed (resp. open) subset of $X$, then $\Phi _{|W}=(\Phi _{|O})_{|C}$ and therefore $\Phi _{|W}$ is a normal and constructible family of supports on $W$. By the above natural isomorphism and Proposition \ref{prop phi soft} it follows that ${\mathcal G}^X$ is $\Phi $-soft. Similarly, we obtain the other implication.

%(ii) Suppose that ${\mathcal H}$ is $\Phi $-soft. Then clearly ${\mathcal H}_{|W}$ is $\Phi _{|W}$-soft and since ${\mathcal H}_W=({\mathcal H}_{|W})^X$ the result follows from (i).
%\qed

A special and useful case of Proposition \ref{prop soft loc and phi-dim} is when $X$ is an object  in the category of o-minimal spectra of definable spaces and  $\Phi $ is a normal and constructible family of supports on $X$.

\begin{prop}\label{prop injective for soft and flabby}
Assume that $X$ is an object  in the category of o-minimal spectra of definable spaces, $\Phi $ is a normal and constructible family of supports on $X$ and $Y$ is a subspace of $X$ such that $D\cap Y$ is a quasi-compact subset for every $D\in \Phi $. Then the
full additive subcategory of ${\rm Mod}(k_Y)$ of   $\Phi \cap Y$-soft $k$-sheaves
is $\Gamma _{\Phi \cap Y}(Y;\bullet )$-injective, i.e.:
\begin{enumerate}
\item
For every ${\mathcal F}\in {\rm Mod}(k_Y)$ there exists a  $\Phi \cap Y$-soft
${\mathcal F}'\in {\rm Mod}(k_Y)$ and an exact sequence
$0\rightarrow {\mathcal F}\rightarrow {\mathcal F}'$.
\item
If $0\rightarrow {\mathcal F}'\rightarrow {\mathcal F}\rightarrow {\mathcal F}''
\rightarrow 0$ is an exact sequence in ${\rm Mod}(k_Y)$ and ${\mathcal F}'$ is
%, ${\mathcal F}$ and ${\mathcal F}''$
 $\Phi \cap Y$-soft, then
$0\into \Gamma  _{\Phi \cap Y}(Y;{\mathcal F}')\into \Gamma  _{\Phi \cap Y}(Y;{\mathcal F})\into
\Gamma  _{\Phi \cap Y}(Y;{\mathcal F}'')\into 0$ is an exact sequence.
\item
If $0\rightarrow {\mathcal F}'\rightarrow {\mathcal F}\rightarrow {\mathcal F}''
\rightarrow 0$ is an exact sequence in ${\rm Mod}(k_Y)$ and ${\mathcal F}'$ and
${\mathcal F}$  are $\Phi \cap Y$-soft, then ${\mathcal F}''$ is  $\Phi \cap Y$-soft.
\end{enumerate}
%Moreover, if $\Phi $ is a normal and constructible  family of supports on $X$ then the same holds with injective or flabby replaced by $\Phi \cap Y$-soft.
\end{prop}

\pf
The result for the full additive subcategory of ${\rm Mod}(k_Y)$ of injective (and flabby) $k$-sheaves  is classical for topological spaces (see for example \cite{ks1}, Proposition 2.4.3).
%Indeed, since there are enough injectives, (1) holds for the injective case; since every injective sheaf is flabby (\cite{b} Chapter II, Proposition 5.3), (1) also holds for the flabby case; (2) and (3) for the flabby (and hence for the injective case) are proved in \cite{b} Chapter II, Theorem 5.4.
Thus (1) holds for the $\Phi \cap Y$-soft case since injective $k$-sheaves are $\Phi \cap Y$-soft.
%On the other hand,  (2)  follows from  the shrinking lemma (\cite{ejp} Proposition 2.17) as in the real algebraic case (\cite{D3} Chapter II, Theorem 4.12, Case 1).
%topological case in \cite{b} Chapter II, Theorem 9.9 to reduce the proof to the case where $X$ is normal and $\Phi $ is the family of all closed subsets of $X$. But then we can finish the proof of (2) by applying \cite{D2} Lemma 2.3 which, as we pointed out in \cite{ejp} page 173, holds in any normal spectral space.
%\cite{ejp} Proposition 3.4 (2).

We now prove (2). Let $s'' \in \Gamma_{\Phi\cap Y}(Y,\F'')$. Then since $\Phi$ is normal and constructible, $\supp\,s'' \subset V$, with $V$ open constructible in $X$ and $\overline{V} \in \Phi$. Now, let us consider the exact sequence
$$
0\to\F'_{V \cap Y}\to\F_{V \cap Y}\to\F''_{V\cap Y}\to 0.
$$
By Proposition \ref{prop soft loc and phi-dim} (iii) we have that $\F'_{Y \cap V}$ is still $\Phi \cap Y$-soft. Replacing $\F',\F,\F''$ with $\F'_{V \cap Y},\F_{V \cap Y},\F''_{V\cap Y}$ we are reduced to prove that the sequence
$$
0\to\Gamma(Y;\F')\to\Gamma(Y;\F)\to\Gamma(Y;\F'')\to 0
$$
is exact when $Y=Y\cap\overline{V}$.
Let $s'' \in \Gamma(Y;\F'')$, and let $\{D_i\}_{i=1}^n$, $D_i \in \Phi \cap Y$ be
a finite covering of $Y$ such that there exists $s_i \in
\Gamma(D_i;\F)$ whose image is $s''|_{D_i}$. There exists such a covering since $\Phi$ is normal and $Y \cap \overline{V}$ is quasi-compact. For $n \geq 2$ on $D_1
\cap D_2$ $s_1-s_2$ defines a section of $\Gamma(D_1 \cap D_2;\F')$
which extends to $s' \in \Gamma(Y;\F')$ since $\F'$ is $\Phi \cap Y$-soft. Replace $s_1$ with
$s_1-s'$. We may suppose that $s_1=s_2$ on $D_1 \cap D_2$. Then
there exists $\widetilde s \in \Gamma(D_1 \cup D_2;\F)$ such that
$\widetilde s|_{D_i}=s_i$, $i=1,2$. Thus the induction proceeds.

Finally, (3) follows at once from (2) by a simple diagram chase using Proposition  \ref{prop phi soft} (3): let $Z$ be a
%constructible
set in $\Phi \cap Y$ and consider the following commutative diagram

\[\xymatrix{
\Gamma _{\Phi \cap Y}(Y;{\mathcal F}) \ar[d]^{\alpha } \ar[rr]^{\delta } &&
\Gamma  _{\Phi \cap Y}(Y;{\mathcal F}'' )\ar[d]^{\gamma }\\
\Gamma  _{\Phi \cap Y\cap Z}(Z;{\mathcal F}_{|Z})\ar[rr]^{\beta }&&
\Gamma  _{\Phi \cap Y\cap Z}(Z;{\mathcal F}''_{|Z})
.}
\]
By hypothesis on ${\mathcal F}$, $\alpha $ is surjective. By (2)
%Since restriction of sheaves on closed subsets is an exact functor, $Z$ is normal and ${\mathcal F}'_{|Z}$ and ${\mathcal F}_{|Z}$ are soft, by  \cite{D2} Lemma 2.4 (which, as we pointed out in \cite{ejp} page 173, holds in any normal spectral space),
%\cite{ejp} Proposition 3.4 (3),
$\beta $ is surjective. Therefore, $\gamma $ is surjective as required.
\qed

Hence, if  $X$ is an object in the category of o-minimal spectra of definable spaces,  $\Phi $ is a normal and constructible family of supports on $X$ and $Y$ is a subspace of $X$ such that $D\cap Y$ is a  quasi-compact  subset for every $D\in \Phi $. Then one can take a $\Phi \cap Y$-soft resolution of ${\mathcal F}$ to compute $H^*_{\Phi \cap Y}(Y;{\mathcal F})$.\\

\begin{expl}\label{expl injective for soft and flabby}
{\em
Some particular cases of Proposition \ref{prop injective for soft and flabby} are:
\begin{itemize}
\item
if $Y=U$ is open constructible such that $\bar{U}\in \Phi$, then the family of $\Phi \cap U$-soft sheaves in $\mod (k_U)$ is $\Gamma (U;\bullet )$-injective.

\item
If $Y=D\in \Phi$, then the family of $\Phi _{|D}$-soft sheaves in $\mod (k_D)$ is $\Gamma (D;\bullet )$-injective.

\end{itemize}
}
\end{expl}

\begin{cor}\label{cor hphi and !}
Assume that $X$ is an object  in the category of o-minimal spectra of definable spaces. Suppose either that   $\Phi $ is a normal and  constructible family of supports on $X$ and  $W$ is a (constructible) locally closed subset of $X$ or that $\Phi$ is any family of supports on $X$ and $W$ is a closed subset of $X$. If $\F\in \mod (k_W)$, then
$$H_{\Phi }^*(X;i_{W!}\F)=H_{\Phi |W}^*(W;\F).$$
\end{cor}

\pf
The second case is covered by \cite{b} Chapter II, 10.1. If $W$ is closed in an open subset $U$ of $X$, then $\Phi _{|U}$ is a normal and constructible family of supports on $U$ and $\Phi _{|W}=\Phi _{|U}\cap W.$ And the result follows from Propositions \ref{prop injective for soft and flabby}, \ref{prop soft loc and phi-dim} (ii) and  \ref{prop section on locally closed}.
\qed

The following will be useful in the next subsection:

\begin{prop}\label{prop soft open and phi-dim}
Assume that  $X$ is an object in the category of o-minimal spectra of definable spaces, ${\mathcal F}$ is a sheaf in ${\rm Mod}(k_X)$ and  $\Phi$ is a normal and constructible family of supports on $X$. The following are equivalent:
\begin{enumerate}
\item
${\mathcal F}$ is $\Phi $-soft;
\item
${\mathcal F}_U$ is $\Gamma _{\Phi }$-acyclic for all open and constructible $U\subseteq X$;
\item
$H^1_{\Phi }(X;{\mathcal F}_U)=0$ for all open and constructible $U\subseteq X$;
%\item
%$H^1_{\Phi |U}(U,{\mathcal F}_{|U})=0$ for all open and constructible $U\subseteq X$.
\end{enumerate}
\end{prop}

\pf
(1) $\Rightarrow $ (2) follows from Propositions \ref{prop injective for soft and flabby} and   \ref{prop soft loc and phi-dim} (iii). (2) $\Rightarrow $ (3) is trivial. To show that (3) implies  (1), consider a constructible closed set $C$ in $\Phi $ and the exact sequence $0\into {\mathcal F}_{X\setminus C}\into {\mathcal F}\into {\mathcal F}_C\into 0$.
% of sheaves in ${\rm Sh}(X)$ (\cite{b} Chapter I, 2.6).
The associated long exact cohomology sequence
$$\dots \rightarrow \Gamma _{\Phi }(X;{\mathcal F})\rightarrow \Gamma _{\Phi }(X;{\mathcal F}_C)\rightarrow H_{\Phi }^1(X;{\mathcal F}_{X\setminus C})\rightarrow \dots $$
shows that $\Gamma _{\Phi }(X;{\mathcal F})\into  \Gamma _{\Phi }(X;{\mathcal F}_C)$ is surjective. Hence ${\mathcal F}$ is $\Phi $-soft by Proposition \ref{prop phi soft} (5).
\qed

\end{subsection}

\begin{subsection}{Cohomological $\Phi $-dimension}\label{subsection cohomo phi-dim}

Recall that for  a topological space  $X$  and $\Phi $ a family of supports on $X$, the {\it cohomological $\Phi $-dimension of $X$} is the smallest $n$ such that $H_{\Phi }^q(X;{\mathcal F})=0$ for all $q>n$ and all sheaves ${\mathcal F}$ in ${\rm Mod}(k_X)$.

%For an object $X$ in the category of o-minimal spectra of definable sets and a constructible family of supports $\Phi $ on $X$, let  $$\dim \Phi=\sup \{\dim Z:Z\in \Phi \,\,{\rm and}\,\, Z\,\,{\rm is}\,\,{\rm constructible}\}.$$ Clearly $\dim \Phi \leq \dim X.$

The following holds:

\begin{prop}\label{clm soft open and phi-dim}
Assume that  $X$ is an object in the category of o-minimal spectra of definable spaces and  $\Phi$ is a normal and constructible family of supports on $X$. Let  ${\mathcal F}$  be a sheaf in ${\rm Mod}(k_X)$. Then the following are equivalent:
\begin{enumerate}
\item
If $0\into {\mathcal F}\into {\mathcal I}^0\into {\mathcal I}^1\into \cdots \into {\mathcal I}^n\into 0$
is an exact sequence of sheaves in ${\rm Mod}(k_X)$ such that ${\mathcal I}^k$ is $\Phi $-soft  for $0\leq k\leq n-1$. Then ${\mathcal I}^n$ is $\Phi $-soft.
\item
${\mathcal F}$ has a $\Phi $-soft resolution of length $n$;
\item
$H^k_{\Phi }(X;{\mathcal F}_U)=H^k_{\Phi |U}(U;{\mathcal F}_{|U})=0$ for all open and constructible $U\subseteq X$ and all $k>n$.
%\item
%$H^{n+1}_{\Phi }(X;{\mathcal F}_{U})=0$ for all open and constructible $U\subseteq X$.
\end{enumerate}
\end{prop}

\pf
%The implications (1) $\Rightarrow $ (2) $\Rightarrow $ (3) $\Rightarrow $ (4) are clear. On the other hand, the implication  (4) $\Rightarrow $ (1)
The result follows from Proposition \ref{prop soft open and phi-dim} (2) and is a particular case of a general result of homological algebra (\cite{ks1}, Exercise I.19): let $F$ be a left exact functor and let ${\mathcal J}$ be the family of $F$-acyclic objects. Suppose that ${\mathcal J}$ is cogenerating. Then (1) $\Leftrightarrow$ (2) $\Leftrightarrow$ (3) with ${\mathcal J}$ instead of $\Phi$-soft and $F$ instead of $\Gamma_\Phi(X;(\bullet )_U)$.
%as in \cite{b} Chapter II, Theorem 16.2.
\qed

\begin{thm}\label{thm cohomo phi dim}
Let  $X$ be an object in the category of o-minimal spectra of definable spaces and let $\Phi $ be a normal and constructible family of supports on $X$. Then the cohomological $\Phi $-dimension of $X$ is bounded by  $\dim X$.
\end{thm}

\pf
To prove our theorem we will use (1) of Proposition \ref{clm soft open and phi-dim}. Let $n=\dim X$. Then, in this situation it suffices to prove that ${\mathcal I}^n_{|Z}$ is soft for every constructible subset $Z$ of $X$ which is in $\Phi $ (Proposition  \ref{prop phi soft} (4)).  Since $\Phi $ is normal, there is a constructible neighborhood $Y$ of $Z$ in $X$ which is in $\Phi $. If we show that ${\mathcal I}^n_{|Y}$ is soft, then it will follow that ${\mathcal I}^n_{|Z}$ is soft (Proposition  \ref{prop phi soft} (2)).

Let $U$ be an open and constructible subset of $Y$. By hypothesis and Proposition \ref{prop soft open and phi-dim} each $({\mathcal I}^k_{|Y})_U$ is acyclic for $0\leq k \leq n-1$. Let ${\mathcal Z}^k={\rm ker}(({\mathcal I}^k_{|Y})_U\into ({\mathcal I}^{k+1}_{|Y})_U).$ Then the long exact cohomology sequences of the short exact sequences $0\into {\mathcal Z}^k\into ({\mathcal I}^k_{|Y})_U\into {\mathcal Z}^{k+1}\into 0$ show that
{\tiny
$$H^q(Y;({\mathcal I}^n_{|Y})_U)=H^q(Y;{\mathcal Z}^n)=H^{q+1}(Y;{\mathcal Z}^{n-1})=\cdots =H^{q+n}(Y;{\mathcal Z}^0)=H^{q+n}(Y;({\mathcal F}_{|Y})_U).$$
}

Since $Y$ is normal, constructible and $\dim Y= n$ we have $H^q(Y;{\mathcal G})=0$ for $q>n$ and every sheaf ${\mathcal G}$ on $Y$ (\cite{ejp} Proposition 4.2). Thus  $H^1(Y;({\mathcal I}^n_{|Y})_U)=0$. Since $U$ was an arbitrary  open and constructible subset of $Y$, it follows from Proposition \ref{prop soft open and phi-dim} that  ${\mathcal I}^n_{|Y}$ is soft as required.
\qed

\begin{prop}\label{prop soft and tensor}
Assume that $X$ is an object  in the category of o-minimal spectra of definable spaces and  $\Phi $ is a normal and  constructible family of supports on $X$. If  ${\mathcal G}\in {\rm Mod}(k_X)$ is  $\Phi$-soft, then  for every   ${\mathcal F}\in {\rm Mod}(k_X)$ we have that  ${\mathcal G}\otimes {\mathcal F}\in {\rm Mod}(k_X)$ is $\Phi $-soft.
\end{prop}

\pf
By Theorem \ref{thm cohomo phi dim}, $X$ has finite cohomological $\Phi$-dimension. Suppose that the cohomological $\Phi$-dimension of $X$ is $n$.
Since the family of the constant sheaves $\{k_U\}$, $U$ constructible open subset of $X$ is generating, there is a resolution of $\F$
$${\mathcal P} _{n-1}\stackrel{\partial _{n-1}}\rightarrow {\mathcal P}_{n-2}\cdots {\mathcal P} _1\stackrel{\partial _1}\rightarrow {\mathcal P}_0\rightarrow \F \rightarrow 0$$
where the ${\mathcal P}_i$'s are direct sums of sheaves of the form $k_U$, $U$ constructible (see \cite{ks1}, Proposition 2.4.12). From Proposition \ref{prop soft loc and phi-dim} (iii) it follows that $\G_U\simeq \G \otimes  k_U$ is $\Phi$-soft. Since the direct sum of $\Phi$-soft sheaves in ${\rm Mod}(k_X)$ is $\Phi$-soft (Corollary \ref{cor lim soft}) each $\G\otimes {\mathcal P}_i $ is $\Phi$-soft.

From the resolution above we obtain an exact sequence of sheaves {\tiny
$$
%0\rightarrow \G\otimes {\mathcal K}\rightarrow
\G\otimes {\mathcal P} _{n-1}\stackrel{\partial _{n-1}}\rightarrow  \G\otimes {\mathcal P}_{n-2}\cdots
 \rightarrow \G\otimes {\mathcal P} _1\stackrel{\partial _1}\rightarrow \G\otimes {\mathcal P}_0\rightarrow  \G\otimes \F \rightarrow 0.$$ }
%Using Theorem \ref{thm cohomo phi dim} and Proposition \ref{clm soft open and phi-dim} we conclude that $\G\otimes {\mathcal F}$ is $\Phi$-soft.
By Proposition \ref{clm soft open and phi-dim}, since $\G \otimes  \mathcal{P}_i$ is $\Phi$-soft for $i=0,\dots,n-2$,
we conclude that $\G\otimes {\mathcal F}$ is $\Phi$-soft.
\qed

\end{subsection}
\end{section}

\begin{section}{Duality with coefficient in a field}\label{section duality with coeff in a field}

In this section we will work in the category of definable spaces with continuous definable maps and $k$-sheaves on such spaces will be considered always relative to the o-minimal site.  In our results we will have a definably normal, definably locally compact definable space $X$ and the family of definable supports $c$ on $X$ of definably compact definable subsets of $X$. By Example \ref{expl supp c} and Remark \ref{nrmk def supp and cons supp} the corresponding constructible family of supports on the o-minimal spectra of $X$ will be a normal and constructible family of supports. Hence, by the tilde isomorphism in the category of $k$-sheaves given by Theorem \ref{thm main iso on sheaves} and our Definition \ref{defn  sheaf cohomo Phi}, in our proofs we will apply the results of Section \ref{section appendix} since they transfer to this  definable setting.

\begin{nrmk}\label{nrmk pv dual for phi}
{\em
We observe that since all the results of this section depend only on Section \ref{section appendix}, they hold on an arbitrary definable space $X$ replacing $c$ by a definably normal family of definable supports $\Phi $ on $X$. In particular, these results hold on any definable space $X$ on which $c$  is a definably normal family of definable supports.
}
\end{nrmk}

\begin{subsection}{Sheaves of linear forms}\label{subsection sheaves of linear forms}

Here we shall work with a fixed field $k$. For a $k$-vector space $N$ we let $N^{\,\vee}$ denote the dual $k$-vector space, i.e. $N^{\, \vee}={\rm Hom}_k(N,k).$

Let $X$ be a definably normal, definably locally compact definable space and ${\mathcal F}$ a $k$-sheaf on $X$. From now on, given a locally closed subset $Z$ of $X$, we will write $\Gamma_c(Z;\F)$ instead of $\Gamma_{c|Z}(Z;\F)$ for short. The inclusion $V\into U$ of two open definable subsets of $X$ will induce a map
$$
\xymatrix{
\Gamma_c(X;\mathcal{F}_V) \ar[r] \ar[d]^\wr & \Gamma_c(X;\mathcal{F}_U) \ar[d]^\wr \\
\Gamma_c(V;\mathcal{F}) \ar[r] & \Gamma_c(U;\mathcal{F})
}
$$
%$$\Gamma _c(V,{\mathcal F})\into \Gamma _c(U,{\mathcal F})$$
``extension by zero''. (Where the vertical isomorphisms are a consequence of Proposition \ref{prop section on locally closed} with $\Phi=c$). The $k$-linear dual of this
$$\Gamma _c(U;{\mathcal F})^{\,\vee}\into \Gamma _c(V;{\mathcal F})^{\,\vee}$$
gives rise to restriction maps in a presheaf ${\mathcal F}^{\, \vee}$ defined by
$$\Gamma (U;{\mathcal F}^{\,\vee})=\Gamma _c(U;{\mathcal F})^{\, \vee}.$$

\begin{prop}\label{prop dual sheaf}
Let $X$ be a definably normal, definably locally compact definable space. For every $c$-soft  $k$-sheaf ${\mathcal F}$ on $X$,  the presheaf ${\mathcal F}^{\,\vee}$ is a sheaf.
\end{prop}

\pf
By Proposition \ref{prop psh and sh}, it is enough to show that for any two open definable subsets $W$ and $V$ of $X$ the sequence
$$0\into \Gamma (V\cup W;{\mathcal F}^{\, \vee})\into  \Gamma (V;{\mathcal F}^{\, \vee})\oplus \Gamma (W;{\mathcal F}^{\, \vee})\into \Gamma (V\cap W;{\mathcal F}^{\, \vee})$$
formed by the sum and difference between two restriction maps is exact.

Consider the Mayer-Vietoris sequence
{\tiny
$$0\into \Gamma _c(V\cap W;{\mathcal F})\into \Gamma _c(V;{\mathcal F})\oplus \Gamma _c(W;{\mathcal F})\into \Gamma _c(V\cup W;{\mathcal F})\into  H_c^1(V\cap W;{\mathcal F})$$
}
and notice that $H_c^1(V\cap W;{\mathcal F})=0$ since the restriction of ${\mathcal F}$ to $V\cap W$ is $c$-soft by Proposition \ref{prop soft loc and phi-dim} (i). The result now follows by taking the $k$-linear dual of the Mayer-Vietoris sequence.
\qed

\begin{prop}\label{prop soft and iso dual}
Let $X$ be a definably normal, definably locally compact definable space. Let ${\mathcal G}$ be a $c$-soft $k$-sheaf on $X$. There is a natural isomorphism
$$\Gamma _c(X;{\mathcal F}\otimes  \G)^{\, \vee}\simeq {\rm Hom}(\F, \G ^{\, \vee})$$
as ${\mathcal F}$ varies through the category of $k$-sheaves on $X$.
\end{prop}

\pf
Let $U$ be an open definable subset. Consider the natural maps
$$\Gamma (U;\F)\otimes \Gamma _c(U;\G)\into \Gamma _c(U;\F\otimes \G )\rightarrow \Gamma _c(X;\F\otimes  \G )$$
The dual of the composite can be written
$$\Gamma _c(X;\F \otimes  \G)^{\, \vee}\into {\rm Hom}(\Gamma (U;\F), \Gamma _c(U;\G)^{\, \vee})$$
By variation of $U$ this defines a map
\begin{equation}\label{eq soft and iso dual}
\Gamma _c(X;\F\otimes  \G)^{\, \vee}\into {\rm Hom}(\F, \G ^{\, \vee})
\end{equation}
which we must show that it is an isomorphism.

(i) First we consider the case where $\F =k_U$ where $U$ is an open definable subset. We have
$$\Gamma _c(X;\G_U)^{\, \vee }= \Gamma _c(U;\G)^{\, \vee }=\Gamma (U;\G ^{\,\vee})={\rm Hom}(k_U, \G ^{\,\vee}).$$
These identifications transform the map (\ref{eq soft and iso dual}) into the identity.

(ii) For the general case, consider a presentation of $\F$ of the form
$${\mathcal P}\into {\mathcal Q}\into \F\into 0$$
where ${\mathcal P}$ and ${\mathcal Q}$ are direct sums of sheaves of the form $k_U$ as above (see \cite{ks1}, Proposition 2.4.12). Let us consider the following diagram with exact rows
{\tiny
\[\xymatrix{
0 \ar[r] & \Gamma _c(X;{\mathcal F}\otimes \G)^{\,\vee } \ar[d]^{ } \ar[r]^{ } & \Gamma _c(X;{\mathcal Q}\otimes {\mathcal G})^{\,\vee } \ar[d]^{} \ar[r]^{} &\Gamma _c(X;{\mathcal P}\otimes {\mathcal G})^{\,\vee } \ar[d]^{}\\
0\ar[r]^{} & {\rm Hom}({\mathcal F},\G ^{\,\vee}) \ar[r]^{ } &   {\rm Hom}({\mathcal Q},\G ^{\,\vee}) \ar[r]^{ } &  {\rm Hom}({\mathcal P},\G ^{\,\vee}) .}
\]}
The two functors of (\ref{eq soft and iso dual}) transform direct sums into direct products. It follows that the two vertical maps to the right are isomorphisms.
Then it follows from the five lemma that the first vertical arrow is an isomorphism.
\qed

\begin{cor}\label{cor soft and iso dual}
Let $X$ be a definably normal, definably locally compact definable space. Let ${\mathcal G}$ be a $c$-soft $k$-sheaf on $X$. Then the sheaf $\G ^{\,\vee}$ is injective in  the category of $k$-sheaves on $X$.
\end{cor}

\pf
By Proposition \ref{prop soft and iso dual}, we must show that $$\F \mapsto \Gamma _c(X;\F\otimes \G)^{\,\vee}$$ is an exact functor. But this follows from Propositions \ref{prop soft and tensor} and \ref{prop injective for soft and flabby} and the exactness of $\vee$ in $\mathrm{Mod}(k)$.
\qed

\end{subsection}

\begin{subsection}{Verdier duality}\label{subsection verdier duality}

If  $X$ is a definably normal, definably locally compact definable space  we will let ${\rm D}^+(k_X)$ denote the derived category of bounded below complexes of $k$-sheaves on $X$. We are now ready to prove our main result:

%In the following we will write $a_{X!}\F$ instead of $\F$.

\begin{thm}[Verdier duality]\label{thm vd dual}
Let $X$ denote a definably normal, definably locally compact definable space. Then there exists an object ${\mathcal D}^*$ in ${\rm D}^+(k_X)$ and a natural isomorphism
$$
{\rm RHom}(\F^*,{\mathcal D}^*) \simeq {\rm RHom}(R\Gamma_c(X;\F^*),k)
$$
as   $\F^*$ varies through  ${\rm D}^+(k_X).$
\end{thm}

\pf
For a complex $L^*$ of $k$-vector spaces we put $L^{*\,\vee}={\rm Hom}^*(L^*,k)$ with the notation of \cite{i} I.4.3. Notice also that $L^{*\,\vee}$ is a complex of $k$-vector spaces  whose $p$'th differential is given by $$(-1)^{p+1}(\partial ^{-p-1})^{\,\vee}:(L^{-p})^{\,\vee }\into (L^{-p-1})^{\, \vee }.$$
This formula will also be used to extend the functor $\G \mapsto \G^{\, \vee }$ on the category of $k$-sheaves on $X$ given by Proposition \ref{prop dual sheaf} to complexes of $k$-sheaves.

By Theorem \ref{thm cohomo phi dim} $X$ has finite cohomological $c$-dimension, hence by Proposition \ref{clm soft open and phi-dim} (1) the constant sheaf $k_X$ admits a bounded $c$-soft resolution $\G^*$.
 %(Which exists since $X$ has finite cohomological $c$-dimension - Theorem \ref{thm cohomo phi dim}).
 By Corollary \ref{cor soft and iso dual}, $\G ^{*\,\vee}$ is a bounded complex of injective $k$-sheaves. For an injective complex $\I ^*$ quasi-isomorphic to $\F^*$ in the derived category of bounded below complexes of $k$-sheaves on $X$  and integers $p$ and $q$ we have, by Proposition \ref{prop soft and iso dual},  a canonical isomorphism
$$\Gamma _c(X;\I ^p\otimes \G ^q)^{\,\vee}={\rm Hom}(\I ^p,\G ^{q\,\vee}).$$
%Taking the direct sum over all $p$ and $q$ such that $p+q=-n$ we find a canonical isomorphism of $k$-vector spaces
%$$[\Gamma _c(X,\I ^*\otimes  \G^*)^{\,\vee}]^n={\rm Hom}^n(\I^*,\G^{*\,\vee})$$
%with the right sign conventions,
giving an isomorphism of complexes
\begin{equation}\label{eq vd dual}
\Gamma _c(X;\I ^*\otimes  \G^*)^{\,\vee}={\rm Hom}^*(\I^*,\G^{*\,\vee}).
\end{equation}
From the quasi-isomorphism $k_X\into \G^*$ we deduce a quasi-isomorpism
%$$\I ^*\into \I^*\otimes \G^*.$$
%Since both complexes are $c$-soft (Proposition %\ref{prop soft and tensor}), we have still another quasi-isomorphism (\cite{i} I.7.5)
%$$\Gamma _c(X,\I^*)\into \Gamma _c(X,\I^*\otimes \G^*)$$
%whose $k$-linear dual
$$\Gamma _c(X;\I^*\otimes \G^*)^{\,\vee}\into \Gamma _c(X;\I^*)^{\,\vee}$$
which yields a final quasi-isomorphism
\begin{equation}\label{eq vd dual dual}
{\rm Hom}^*(\I^*,\G^{*\,\vee})\into \Gamma _c(X;\I^*)^{\,\vee}.
\end{equation}
%Passing to $H^0$ we obtain
%$$[\I^*,\G^{*\,\vee}]\simeq H^0\Gamma _c(X,\I^*)^{\,\vee}\simeq [\Gamma _c(X,\I^*),k].$$
Finally put ${\mathcal D}^*=\G^{*\,\vee}.$
\qed

The complex ${\mathcal D}^*$ above is called the {\it dualizing complex}. It is  a bounded below complex of injective $k$-sheaves uniquely  determined up to homotopy  and so the cohomology $k$-sheaves ${\mathcal H}^p{\mathcal D}^*,\,\,p\in {\ZZ}$, are uniquely determined up to isomorphism.

The proof  above shows the following:

\begin{nrmk}\label{nrmk dual complex}
{\em
The dualizing complex for a definably normal, definably locally compact definable space of cohomological $c$-dimension $n$ can be represented by a complex ${\mathcal D}^*$ of injective $k$-sheaves where
$${\mathcal D}^i=0\,\,\,\,\,{\rm for}\,\,\,\,\,i\notin [-n,0].$$
}
\end{nrmk}

Recall that the inclusion $V\into U$ of open definable subsets of $X$ give rise to the extension by zero map
$$H_c^p(V;k_X)\into H_c^p(U;k_X)$$
whose $k$-linear dual
$$H_c^p(U;k_X)^{\,\vee}\into H_c^p(V;k_X)^{\,\vee}$$
gives rise to a presheaf $U\mapsto H_c^p(U;k_X)^{\,\vee}.$

\begin{prop}\label{prop coho dual complex}
Let ${\mathcal D}^*$ denote the dualizing complex for the definably normal, definably locally compact definable space $X$. For any integer $p$, the cohomology $k$-sheaf ${\mathcal H}^{-p}{\mathcal D}^*$ is the sheaf associated to the $k$-presheaf
$$U\mapsto H_c^p(U;k_X)^{\,\vee}.$$
\end{prop}

\pf
Recall the isomorphism $H_c^p(U;k_X) \simeq H_c^p(X;k_U)$. Passing to the dual and using Theorem \ref{thm vd dual} we have the chain of isomorphisms
$$
H^p_c(U;k_X)^\vee \simeq H_c^p(X;k_U)^\vee \simeq H^{-p}{\rm Hom}(k_U,{\mathcal D}^*) \simeq H^{-p}(U;{\mathcal D}^*)
$$
and the result follows since ${\mathcal H}^{-p}{\mathcal D}^*$ is the $k$-sheaf associated to the $k$-presheaf $U \mapsto H^{-p}(U;{\mathcal D}^*)$.
\qed

%First we observe that there is a canonical isomorphism for every $p\in {\ZZ}$
%\begin{equation}\label{eq coho dual complex1}
%H_c^p(X,{\mathcal F})^{\,\vee}=H^{-p}{\rm Hom}({\mathcal F},{\mathcal D}^*)
%\end{equation}
%as ${\mathcal F}$ runs through ${\rm Sh}_{{\rm dtop}}(X,k).$ For this, consider an injective resolution $\I^*$ of $\F.$ Then
%$$H_c^p(X,\F)=H^{-p}\Gamma _c(X,\I^*)=H^0\Gamma _c(X,\I^*[p])$$
%and from the duality formula (Theorem \ref{thm %vd dual}) we get
%$$H_c^p(X,\F)^{\,\vee}=[\Gamma _c(X,\I^*[p]),k]=[\I^*[p],{\mathcal D}^*]=[\I^*,{\mathcal D}^*[-p]].$$
%Now the quasi-isomorphism $\F\into \I^*$ allows us to make the identification
%$$[\I^*,{\mathcal D}^*[-p]]=H^{-p}{\rm Hom}^*(\I^*,{\mathcal D}^*)=H^{-p}{\rm Hom}(\F,{\mathcal D}^*)$$
%which proves (\ref{eq coho dual complex1}).

%Now take the special case $\F=j_!\underline{k}$ where $j:U\into X$ is the inclusion of an open definable subset:
%$$H_c^p(U,\underline{k})^{\,\vee}=H_c^p(X,j_!\underline{k})^{\,\vee}=H^{-p}{\rm Hom}(j_!\underline{k},{\mathcal D}^*).$$
%Using the adjunction formula and \cite{i} II.7.2 (in $\tilde{X}$)
%$${\rm Hom}(j_!\underline{k}, {\mathcal D}^*)={\rm Hom}(\underline{k},j^*{\mathcal D}^*)=\Gamma (U,{\mathcal D}^*)$$
%we get the formula
%\begin{equation}\label{eq coho dual complex2}
%H^{-p}\Gamma (U,{\mathcal D}^*)\simeq H_c^p(U,\underline{k})^{\,\vee}.
%\end{equation}
%Since (\ref{eq coho dual complex1}) is functorial in $\F$, (\ref{eq coho dual complex2}) is compatible with restriction to an open definable subset of $U$.
%\qed

\begin{cor}\label{cor coho dual complex}
On a definably normal, definably locally compact definable space $X$ of cohomological $c$-dimension $n$, the $k$-presheaf
$$U\mapsto H_c^n(U;k_X)^{\,\vee}$$
is a $k$-sheaf.
\end{cor}

\pf
By Remark \ref{nrmk dual complex} we have an exact sequence
$$0\into \Gamma (U;{\mathcal H}^{-n}{\mathcal D}^*)\into \Gamma (U;{\mathcal D}^{-n})\into \Gamma (U;{\mathcal D}^{-n+1}).$$
On the other hand $H^{-n}(U;\mathcal{D}^*)=\ker(\Gamma(U;\mathcal{D}^{-n}) \to \Gamma(U;\mathcal{D}^{-n+1})$.
Moreover, as we saw above
$$H^{-n}(U;{\mathcal D}^*) \simeq H_c^n(U;k_X)^{\,\vee}.$$
Then $\Gamma(U;\mathcal{H}^{-n}\mathcal{D}^*)\simeq H_c^n(U;k_X)^{\,\vee}$ and the result follows.
%Combining these two results, we can see $H_c^n(U,k_X)^{\,\vee}$ as sections in the $k$-sheaf ${\mathcal H}^{-n}{\mathcal D}^*.$ (The $k$-presheaf $U\mapsto H_c^n(U,k_X)^{\,\vee}$ is the kernel of the morphism ${\mathcal D}^{-n}\into {\mathcal D}^{-n+1}$ of $k$-sheaves, so it is a $k$-sheaf.)
\qed

\end{subsection}

\begin{subsection}{Poincar\'e and Alexander duality}\label{subsection poincare duality}

Here we derive Poincar\'e and Alexander duality from the Verdier duality.

\begin{defn}\label{defn orientable def man}
{\em
Let $X$ be a definably normal,  definably locally compact definable manifold of dimension $n$. We say that $X$ {\it has an orientation $k$-sheaf} if for every open definable subset $U$ of $X$ there exists a finite cover of $U$ by open definable subsets $U_1,\dots , U_\ell$ of $U$ such that for each $i$ we have
\begin{equation*}
H_c^p(U_i;k_X)=
\begin{cases}
k \qquad \textmd{if} \qquad p=n\\
\\
\,\,\,\,\,\,\,\,\,\,\,\,\,\,\,\,\,\,\,\,\,\,\,\,\,\,\,\,\,\,\,
\\
0\qquad \textmd{if} \qquad p\neq n.
\end{cases}
\end{equation*}
If $X$ has an orientation sheaf, we call the $k$-sheaf ${\mathcal Or}_X$ on $X$ with sections
$$\Gamma (U;{\mathcal  Or}_X)=H_c^n(U;k_X)^{\, \vee}$$ the {\it orientation $k$-sheaf on $X$}. By Theorem \ref{thm cohomo phi dim}, the cohomological $c$-dimension of $X$ is $\leq n$ and $H_c^n(U_i;k_X) =H_c^n(X;k_{U_i}) \neq 0$ for $i=1,\dots,\ell$, hence  $X$ must have cohomological $c$-dimension $n$. So ${\mathcal  Or}_X$ is indeed a $k$-sheaf on $X$ by Corollary \ref{cor coho dual complex}).

Note also that, since the o-minimal spectra $\tilde{X}$ of $X$ is a quasi-compact (spectral) topological  space, $X$ has an orientation $k$-sheaf if and only if for every $\beta \in \tilde{X}$ and every open definable subset $V$ of $X$ such that $\beta \in \tilde{V}$, there is an open definable subset $U$ of $V$ such that $\beta \in \tilde{U}$ and
\begin{equation*}
H_c^p(U;k_X)=
\begin{cases}
k \qquad \textmd{if} \qquad p=n\\
\\
\,\,\,\,\,\,\,\,\,\,\,\,\,\,\,\,\,\,\,\,\,\,\,\,\,\,\,\,\,\,\,
\\
0\qquad \textmd{if} \qquad p\neq n.
\end{cases}
\end{equation*}
}
\end{defn}

\begin{expl}\label{expl orient sheaf}
{\em
Suppose that ${\mathcal M}$ is an  o-minimal expansion of an ordered field. Let $X$ be a Hausdorff definable manifold of dimension $n$. Since then $X$ is affine and every definable set is definably normal, $X$ is definably normal (\cite{vdd} Chapter 6, Lemma 3.5). Since also $X$ and any open definable subset of $X$ can be covered by finitely many definable sub-balls (\cite{e3} Theorem 1.2),
$X$ is definably locally compact and, computing the o-minimal cohomology with definably compact supports of definable sub-balls, it follows that $X$ has an orientation $k$-sheaf. Observe that the result on coverings by definable sub-balls is related to \cite{bo} Theorem 4.3 (and can be read off from  the proofs of Lemmas 4.1 and 4.2 there) and also to Wilkie's result (\cite{w} Theorem 1.3) which says that every bounded open definable set can be covered by finitely open cells.
}
\end{expl}

Let $X$ be a definably normal,  definably locally compact definable manifold of dimension $n$ with an orientation $k$-sheaf ${\mathcal Or}_X$.  Then the $k$-sheaf ${\mathcal Or}_X$  is locally isomorphic to $k_X$.
%Also, by going to the category of o-minimal spectra and applying \cite{i} II.7.5,

%We have
%$$H^p(X,{\mathcal Or}_X)={\rm Ext}^p(k_X,{\mathcal Or}_X) \qquad ; p\in {\ZZ}.$$

%where ${\rm Ext}$ is calculated in ${\rm Sh}_{{\rm dtop}}(X,k)$. From general principles \cite{i} I.8.6

%We deduce for $p, q\in {\ZZ}$ a cup product
%$$\alpha \cup \beta \in H_c^{p+q}(X,{\mathcal Or}_X) \qquad ; \alpha \in H^p(X,{\mathcal Or}_X), \,\, \beta \in H_c^q(X,k_X)$$
%subject to Poincar\'e duality:

\begin{thm}[Poincar\'e duality]\label{thm poinc dual}
Let $X$ be a definably normal,  definably locally compact definable manifold of dimension $n$ with an orientation $k$-sheaf ${\mathcal Or}_X$. There exists  an isomorphism
$$H^p(X;{\mathcal Or}_X)\into H_c^{n-p}(X;k_X)^{\, \vee}.$$
\end{thm}

\pf
%First note that by our Definition \ref{defn  sheaf cohomo Phi}, the support of a $k$-sheaf on $X$ relative to the o-minimal site is the support of its tilde in the o-minimal spectrum on $X$. Thus, by going to the category of o-minimal spectra of definable spaces and taking stalks, we see that
Proposition \ref{prop coho dual complex} and the fact that $X$ has an orientation $k$-sheaf, imply that
$${\mathcal H}^{-p}{\mathcal D}^*=0\qquad ; p\neq n.$$
On the other hand, by Corollary \ref{cor coho dual complex} we have ${\mathcal H}^{-n}{\mathcal D}^*={\mathcal Or}_X$. Thus we have a quasi-isomorphism
\begin{equation}\label{eq pd or}
{\mathcal Or}_X[n]\simeq {\mathcal D}^*
\end{equation}
%and ${\mathcal Or}_X\into {\mathcal D}^*$ is an injective resolution of the sheaf ${\mathcal Or}_X$.
Therefore we have
$$H^p(X;{\mathcal Or}_X)\simeq H^{p-n} (X;{\mathcal D}^*)  \simeq H^{p-n} {\rm Hom}(k_X,\mathcal{D}^*).$$
By  Verdier duality (Theorem \ref{thm vd dual}) with $\F^*=k_X$ the later is also isomorphic to $H_c^{n-p}(X;k_X)^\vee.$
%$$H^p(X,{\mathcal Or}_X)=H^p\Gamma (X,{\mathcal D}^*[-n])=H^0\Gamma (X,{\mathcal D}^*[p-n])=[k,{\mathcal D}^*[p-n]].$$ Hence, if $k_X\into K^*$ is an injective resolution of $k_X$, we also have
%$$H^p(X,{\mathcal Or}_X)=[{\mathcal K}^*,{\mathcal D}^*[p-n]]=[{\mathcal K}^*[n-p],{\mathcal D}^*].$$
%Applying the characterization of the Verdier duality isomorphism given in Remark \ref{nrmk trace map} (\ref{eq vd trace map dual}) to ${\mathcal I}^*={\mathcal K}^*[n-p]$ we get
%$$[{\mathcal K}^*[n-p],{\mathcal D}^*]=[\Gamma _c(X,{\mathcal K}^*[n-p]),k]=H^{n-p}_c(X,k_X)^{\, \vee}.$$

%Finally we use  the characterization of the Verdier duality isomorphism given in Remark \ref{nrmk trace map} (\ref{eq vd trace map dual}) to identify this isomorphism with the one described by the cup product.
\qed

\begin{defn}\label{defn orientable}
{\em
Let $X$ be a definably normal,  definably locally compact definable manifold of dimension $n$ with an orientation $k$-sheaf ${\mathcal Or}_X$. By a {\it $k$-orientation} we understand an isomorphism
$$k_X\simeq {\mathcal Or}_X$$
of $k$-sheaves. We shall say that $X$ is {\it $k$-orientable}  if a $k$-orientation exists and {\it $k$-unorientable} in the opposite case.
}
\end{defn}

\begin{prop}\label{prop orient unorient}
Let $X$ be a definably connected, definably normal,  definably locally compact definable manifold of dimension $n$ with an orientation $k$-sheaf ${\mathcal Or}_X$. Then
\begin{enumerate}
\item $H^n_c(X;k_X)\simeq k$ if $X$ is $k$-orientable.
\item $H^n_c(X;k_X)\simeq 0$ if $X$ is $k$-unorientable.
\end{enumerate}
\end{prop}

\pf
Since $X$ is definably normal and definably connected, Proposition 4.1 in \cite{ejp} implies that $H^0(X;k_X)=k$ and so (1) follows at once from the Poincar\'e duality (Theorem \ref{thm poinc dual}).

For (2), suppose that $H^n_c(X;k_X)\neq 0$. Then by Theorem \ref{thm poinc dual} there is a non trivial section $s$ of ${\mathcal Or}_X$ over $X$. By our Definition \ref{defn  sheaf cohomo Phi}, the support of $s$ is a closed subset of the o-minimal spectrum of $X$. Since ${\mathcal Or}_X$ is locally isomorphic to $k_X$ it follows that the support of $s$ is also an open subset of the o-minimal spectrum of $X$. But since the o-minimal spectrum of $X$ is connected (Theorem \ref{thm main normal def normal}) it follows that the support of $s$ is the o-minimal spectrum of $X$. Thus ${\mathcal Or}_X\simeq k_X$.
\qed

%Consider the inclusion $i_Z:Z\into X$ of a closed definable subset $Z$ of a definable space $X$.

%By going to the category of o-minimal spectra of definable spaces and applying \cite{i} II.9.8, we have
%$$H^p_Z(X,k_X)={\rm Ext}^p(i_*k_X,{\mathcal F}) \qquad ; p\in {\ZZ}$$

%$$H^p_Z(X,\F)={\rm Ext}^p(k_Z,{\mathcal F}) \qquad ; p\in {\ZZ}$$

%where ${\rm Ext}$ is calculated in ${\rm Sh}_{{\rm dtop}}(X,k)$.

%If $X$ is a definably normal,  definably locally compact definable space,

%then from general principles \cite{i} V.6

%we deduce for $p, q\in {\ZZ}$ a cup product
%$$\alpha \cup \beta \in H_c^{p+q}(X,{\mathcal F}) \qquad ; \alpha \in H^p_Z(X,{\mathcal F}), \,\, \beta \in H_c^q(Z,k_X)$$
%subject to Alexander duality:

\begin{thm}[Alexander duality]\label{thm alex dual}
Let $X$ be a definably normal,  definably locally compact, $k$-orientable  definable manifold of dimension $n$. For $Z$ a closed definable subset of $X$ there exists  an isomorphism
$$H^p_Z(X;k_X)\into H_c^{n-p}(Z;k_X)^{\, \vee}.$$
\end{thm}

\pf
By (\ref{eq pd or}) we have
$H^p_Z(X;k_X)  \simeq H_Z^{p-n}(X;\mathcal{D}^*) \simeq H^{p-n} {\rm Hom}(k_Z,\mathcal{D}^*).$
By  Verdier duality (Theorem \ref{thm vd dual}) with $\F^*=k_Z$ the later is also isomorphic to $H_c^{n-p}(X;k_Z)^{\, \vee} \simeq H^{n-p}_c(Z;k_X)^{\, \vee}.$
%The result follows from  once we observe that the duality isomorphism (\ref{eq vd trace map dual})
%$$\phi \mapsto \int _X \circ \Gamma _c(X,\phi )\qquad ; \phi \in [{\mathcal I}^*,{\mathcal D}^*]\simeq [\Gamma _c(X,{\mathcal I}^*,k]$$
%can be rewritten
%$$[{\mathcal I}^*,{\mathcal D}^*[p-n]]\simeq H^{n-p}\Gamma _c(X,{\mathcal I}^*)^{\, \vee},$$
%and so when applied to an injective resolution $i_*k_X\into I^*$ gives
%$$[i_*k_X,{\mathcal D}^*[p-n]]\simeq H_c^{n-p}(Z,k_X)^{\, \vee}.$$
\qed

\end{subsection}

\begin{subsection}{Duality in o-minimal expansions of fields}\label{subsection duality in fields}

In this subsection  we assume   that the o-minimal structure ${\mathcal M}$ is an  expansion of an ordered field.

Let $X$ be a Hausdorff definable manifold of dimension $n$. Then has we saw in Example \ref{expl orient sheaf} $X$ is affine,  definably normal with an orientation $k$-sheaf.

In o-minimal expansions of fields we have o-minimal singular homology and cohomology theories satisfying the Eilenberg-Steenrod axioms adapted to the o-minimal site (\cite{ew}, \cite{Wo}). By \cite{ew} the o-minimal singular cohomology theory with coefficients in a field $k$ is isomorphic to the o-minimal sheaf cohomology theory with coefficients in the constant sheaf $k_X$. Because of this isomorphism, below we will use the standard notation  from o-minimal singular cohomology and write $k$ for $k_A$
and
$$H^*(A,B;k) \,\,\,\,\,{\rm for}\,\,\,\,\, H^*_{A\setminus B}(A;k)$$
where $B\subseteq A\subseteq X$ are definable subsets of $X$.

O-minimal singular homology theory can be used to obtain an orientation theory for definable manifolds (\cite{bo}, \cite{beo}). (In the papers \cite{bo} and \cite{beo}, orientation is defined by taking homology with coefficients in ${\ZZ}$ but replacing ${\ZZ}$ by $k$ and considering homology groups as $k$-vector spaces one gets the theory of $k$-orientations.) Our goal here is to show an Alexander duality for homology and to conclude that the two orientation theories agree.

First observe that if  $B\subseteq A$ are definably locally closed definable subsets of $X$, then
\begin{equation}\label{eq cohomo comp as limit of closed}
H_c^*(A\setminus B;k)=\varinjlim_{B\subseteq C\subseteq A, \,\, C\,\,{\rm closed},\,\,\bar{A\setminus C}\in c}H^*(A,C; k).
\end{equation}
Let $\Lambda $ be the directed system of definably locally closed subsets $D$ of $A$ such that $B\subseteq D\subseteq A$ and $\bar{A\setminus D}\in c$, directed by reverse inclusion. Since the map that sends $D\in \Lambda $ into $\bar{D}$ is cofinal (even surjective) in the   directed system of definable  closed subsets $C$ of $A$ such that $B\subseteq C\subseteq A$ and $\bar{A\setminus C}\in c$, directed by reverse inclusion, it follows that to prove (\ref{eq cohomo comp as limit of closed}) it is enough to show that
$$H_c^*(A\setminus B;k)=\varinjlim_{D\in \Lambda }H^*(A,D;k),$$
i.e.,   we have to show that the natural homomorphism
$$\varinjlim_{B\subseteq U\subseteq A, \,\,A\setminus U\in c}H^*(A,U;k)\into
\varinjlim_{D\in \Lambda }H^*(A,D;k)
$$
is an isomorphism. But this is a consequence of the following. If $D\in \Lambda $, then there exists an open definable subset $O$ of $A$ such that $D$ is closed in $O$. So, by \cite{vdd} Chapter VIII, 3.3 and 3.4, there is an open definable neighborhood $U$ of $D$ in $O$ such that $D$ is a definable deformation retract of $U$. Therefore, the inclusion $D\into U$ induces an isomorphism $H^*(A, U;k)\into H^*(A,D;k)$.

We are now ready to show the  Alexander duality for o-minimal homology. This is the o-minimal version of \cite{d} Chapter VIII, Theorem 7.14 and the generalization of Theorem 3.5 in \cite{ew2}.

\begin{thm}\label{thm duality for hc}
Let $X$ be a definable manifold of dimension $n$ which is $k$-orientable with respect to homology. Let  $L\subseteq K\subseteq X$ be closed definable sets with $K-L$ closed in $X-L$. Then there is an isomorphism
$$H_c^q(K\setminus L;k)\into H_{n-q}(X\setminus L,X\setminus K;k)$$
for all $q\in {\ZZ}$ which is natural with respect to inclusions.
\end{thm}

\pf
%By  (\ref{eq cohomo comp rel vs non rel}),  $H_c^*(K,L;k)\simeq H_c^*(K\setminus L;k)$.
Let $K'=K\setminus L$, $X'=X\setminus L$, $A$ a definable closed subset of $K'$ such that $\bar{K'\setminus A}\in c$ and $C=\bar{K'\setminus A}$. Then we have the following commutative diagram
{\tiny
\[\xymatrix{
H^{q}(K',A;k)\ar[d] \ar[rr] && H_{n-q}(X'\setminus A,X'\setminus K';k)\ar[d] \\
H^{q}(K'\cap C,A\cap C;k)\ar[rr]^{\cap \zeta '_{K'\cap C}} && H_{n-q}(X'\setminus A\cap C,X'\setminus K'\cap C;k)
.}
\]
}
\noindent
where the vertical arrows are the inclusion homomorphisms which, by the excision
axiom, are isomorphisms. The bottom arrow is the isomorphism of Theorem 3.5 in \cite{ew2}. This diagram goes to the limit to give the isomorphism of the theorem by (\ref{eq cohomo comp as limit of closed}) and
$$H_*(X',X'-K';k)= \varinjlim_{A\subseteq K',\,\,A\,\,{\rm closed},\,\,\bar{K'\setminus A}\in c} H_*(X'-A,X'-K';k)$$
(as $X'=\cup \{X'-A:A\subseteq K',\,\,A\,\,{\rm closed},\,\,\bar{K'\setminus A}\in c \}$).%; compare with \cite{d} Chapter VIII, 5.22).
\qed

Combining Alexander duality for homology (Theorem \ref{thm duality for hc}) and for cohomology  (Theorem \ref{thm alex dual}) we show:

\begin{cor}\label{cor hom cohm orient}
Let $X$ be a Hausdorff definable manifold. Then $X$ is $k$-orientable with respect to homology if and only if $X$ is $k$-orientable with respect to cohomology.
\end{cor}

\pf
Indeed, let $X$ be a Hausdorff definable manifold of dimension $n$. If $X$ is $k$-orientable with respect to homology, then Theorem \ref{thm duality for hc}  implies that for every definably connected, definably compact definable subset $K$ of $X$ we have an isomorphism $H_{n}(X,X\setminus K;k)\simeq k$ which is compatible with inclusions. Applying the dual universal coefficients theorem and going to the limit we obtain $H_c^n(X;k)\simeq k$ showing that $X$ is $k$-orientable (Proposition \ref{prop orient unorient}). If $X$ is $k$-orientable with respect to cohomology, then  Theorem \ref{thm alex dual}  applied to $K$ and  $X$ implies that for every definably connected, definably compact definable subset $K$ of $X$ we have an isomorphism $H^{n}(X,X\setminus K;k)\simeq k$ which is compatible with inclusions. Applying the dual universal coefficients theorem we get an isomorphism $H_{n}(X,X\setminus K;k)\simeq k$ compatible with inclusions which allows us to define a $k$-orientation for $X$ relative to homology.
\qed

\end{subsection}
\end{section}

\end{document}